\documentclass[11pt]{article}
\usepackage{amsfonts}
\usepackage{graphics}
\usepackage{indentfirst}
\usepackage{color}
\usepackage{cite}
\usepackage{latexsym}
\usepackage[paper=a4paper, left=1.6cm, right=1.6cm, top=1.8cm, bottom=1.6cm, headheight=5.5pt, footskip=0.8cm, footnotesep=0.8cm, centering, includefoot]{geometry}
\usepackage{amsmath}
\allowdisplaybreaks
\usepackage{amssymb}
\usepackage[colorlinks, linkcolor=red]{hyperref}
\hypersetup{urlcolor=red, citecolor=red}
\usepackage[dvips]{epsfig}
\usepackage{amscd}

\newtheorem{theorem}{Theorem}[section]
\newtheorem{remark}{Remark}[section]

\newtheorem{lemma}{Lemma}[section]
\newtheorem{corollary}{Corollary}[section]
\newtheorem{proposition}{Proposition}[section]

\DeclareMathOperator{\curl}{curl}

\DeclareMathOperator{\divv}{div}

\makeatletter
\@addtoreset{equation}{section}
\makeatother
\makeatletter
\@addtoreset{equation}{section}
\makeatother

\title{Global existence and decay estimates of strong solutions for compressible non-isentropic magnetohydrodynamic flows with vacuum
\thanks{This research was partially supported by National Natural Science Foundation of China (Nos. 12371227, 11901474, 11901288), Scientific Research Foundation of Jilin Provincial Education Department (No. JJKH20210873KJ), Postdoctoral Science Foundation of China (No. 2021M691219), and Exceptional Young Talents Project of Chongqing Talent (No. cstc2021ycjh-bgzxm0153).
}
}

\author{Yang Liu$\,^{\rm 1}\,$,\ Xin Zhong$\,^{\rm 2}\,${\thanks{Corresponding author. E-mail addresses:
 liuyang0405@ccsfu.edu.cn (Y. Liu),  xzhong1014@amss.ac.cn (X. Zhong).}}
\date{}\\
\footnotesize $^{\rm 1}\,$ College of Mathematics, Changchun Normal
University, Changchun 130032, P. R. China\\
\footnotesize $^{\rm 2}\,$ School of Mathematics and Statistics, Southwest University, Chongqing 400715, P. R. China}

\date{ }

\begin{document}
\maketitle

\begin{abstract}
We study the Cauchy problem of three-dimensional compressible non-isentropic magnetohydrodynamic (MHD) fluids with both interior and far field vacuum states. Applying delicate energy estimates, initial layer analysis, and continuation theory, we establish the global existence and uniqueness of strong solutions, which may be of possibly large oscillations, provided that the initial data are of small total energy. Furthermore, we also derive algebraic decay estimates of the solution. This improves our previous results (Z. Angew. Math. Phys. 71: Paper No. 188, 2020; J. Differential Equations 336: 456--478, 2022) where viscosity coefficients should additionally satisfy $3\mu>\lambda$. To our best knowledge, this is the first result on time-decay rates of solutions to the Cauchy problem of multi-dimensional compressible non-isentropic MHD equations with vacuum.
\end{abstract}



\section{Introduction and main results}

The compressible MHD equations govern the motion of electrically conducting fluids such as plasmas, liquid metals, and electrolytes. They consist of a coupled system of compressible Navier-Stokes equations of fluid dynamics and Maxwell's equations of electromagnetism. Besides their wide physical applicability (see, e.g., \cite{D2017}), the MHD system are also of great interest in mathematics. As a coupled system, the issues of well-posedness and dynamical behaviors of compressible MHD equations are rather complicated to investigate because of the strong coupling and interplay interaction between the fluid motion and the magnetic field.
In this paper, we are concerned with three-dimensional compressible non-isentropic MHD equations which take the following form (see \cite[Chapter 3]{LQ2012}):
\begin{align}\label{a1}
\begin{cases}
 \rho_t+\divv(\rho u)=0,\\
\rho u_t+\rho u\cdot\nabla u-\mu\Delta u-(\lambda+\mu)\nabla\divv u+\nabla p=\curl H\times H,\\
c_v\rho(\theta_t+u\cdot\nabla\theta)
+p\divv u-\kappa\Delta\theta=\lambda(\divv u)^2+\frac{\mu}{2}|\nabla u+(\nabla u)^{tr}|^2+\eta|\curl H|^2,\\
H_t-\eta\Delta H=\curl(u\times H),\\
\divv H=0.
\end{cases}
\end{align}
The unknowns $\rho$, $u=\big(u^1, u^2, u^3\big)$, $H=\big(H^1, H^2, H^3\big)$, and $\theta$ denote the density, the velocity, the magnetic
field, and the temperature, respectively. $p=R\rho\theta$, with positive constant $R$, is the pressure.
The constant viscosity coefficients $\mu$ and $\lambda$ satisfy the physical restrictions
\begin{align}\label{a0}
\mu>0, \ 2\mu+3\lambda\ge 0.
\end{align}
The positive constants $c_v$ and $\kappa$ are the heat capacity and the ratio of the heat conductivity coefficient over the heat capacity, respectively, while $\eta>0$ is the magnetic diffusivity.

We consider the Cauchy problem of \eqref{a1} with the initial condition
\begin{align}\label{a2}
(\rho, u, \theta, H)(x, 0)=(\rho_0, u_0, \theta_0, H_0)(x), \ \ x\in\mathbb{R}^3,
\end{align}
and the far field behavior
\begin{align}\label{a3}
(\rho, u, \theta, H)(x, t)\rightarrow (0, 0, 0, 0)\ \ {\rm as}\ |x|\rightarrow \infty, \ t>0.
\end{align}

When we take no account of the energy equation \eqref{a1}$_3$, the system \eqref{a1} is reduced to the well-known compressible isentropic MHD equations. The mathematical results concerning the global existence of (weak, strong or
classical) solutions to this model can refer for example to \cite{FL2020,SH12,HHPZ,HW10,LXZ13,WW17,LSX16}. In contrast to the isentropic case, the non-isentropic model \eqref{a1} is more in line with reality but the problem  becomes challenging. Hu and Wang \cite{HW08} proved the global existence of weak solutions to the initial boundary value problem of \eqref{a1} with general large initial data and vacuum in terms of Lions-Feireisl compactness framework \cite{F2004,L1998} provided the adiabatic exponent is suitably large (see also \cite{DF2006,LG2014} for related works). Very recently, Li and Sun \cite{LS2021} showed global-in-time weak solutions for the full compressible nonresistive MHD system (that is, $\eta=0$ in \eqref{a1}) in multi-dimensional
spaces with large initial data. Yet the uniqueness and regularity of such weak solutions is still open. Meanwhile, for the initial data satisfying some compatibility condition, Fan and Yu \cite{FY09} established the local existence and uniqueness of strong solutions to the problem \eqref{a1}--\eqref{a3}. Liu and Zhong \cite{LZ20} extended this local existence result to be a global one provided that $\|\rho_0\|_{L^\infty}+\|H_0\|_{L^3}$ is suitably small and viscosity coefficients satisfy $3\mu>\lambda$. However, it eliminates the possibly large oscillations of the initial data since the initial density must be pointwisely small. Later on, Hou-Jiang-Peng \cite{HJP22} showed
global strong solutions under the condition that $\|\rho_0\|_{L^1}+\|H_0\|_{L^2}$ is properly small. Very recently, Liu and Zhong \cite{LZ22} derived global strong solutions provided that some scaling invariant quantity is small enough. In addition, other mathematical topics for compressible non-isentropic MHD equations have also been widely investigated, such as low Mach number limit \cite{JJL12,JJLX14}, blow-up criteria of solutions \cite{HL13,Z19,W21}, and so on. It should be noted that the system \eqref{a1} becomes the compressible non-isentropic Navier-Stokes equations when there is no electromagnetic field, which is one of the most important systems in fluid dynamics, the readers can refer to \cite{HL18,WZ17,JL19,L21,XZ21,YZ20} for studies on the global existence of strong solutions.

In the present paper we further investigate the global existence and uniqueness of strong solutions to the compressible non-isentropic MHD equations in $\mathbb{R}^3$ with vacuum at infinity. Our result reveals that the small initial energy is sufficient to guarantee the global well-posedness for the problem \eqref{a1}--\eqref{a3}. Moreover, we also show that the solutions converge algebraically to zero in $L^2$ as time goes to infinity. It is worth mentioning that the main novelty of this article consists in allowing possibly large oscillations of the initial data as well as obtaining algebraic decay estimates of the solution. To the best of our knowledge, this is the first result on time-decay rates of solutions to the Cauchy problem of multi-dimensional compressible non-isentropic MHD equations with vacuum.

Before formulating our main results, we first explain the notations and conventions used throughout this paper.
 For simplicity, in what follows, we denote
\begin{align*}
\int_{\mathbb{R}^3}\cdot dx=\int \cdot dx.
\end{align*}
For $1\le p\le \infty$ and integer $k\ge 0$, we shall use the standard homogeneous and inhomogeneous
Sobolev spaces as follows:
\begin{align*}
\begin{cases}
L^p=L^p(\mathbb{R}^3),~ ~ W^{k, p}=L^p\cap D^{k, p},~~ H^k=W^{k, 2}, \\
  D^{k, p}=\{u\in L_{loc}^1(\mathbb{R}^3): \|\nabla^ku\|_{L^p}<\infty\}, ~D^k=D^{k, 2},\\
  D_0^1=\{u\in L^6(\mathbb{R}^3): \|\nabla u\|_{L^2}<\infty\}.
\end{cases}
\end{align*}

We can now state our first main result concerning the global existence of strong solutions to the problem \eqref{a1}--\eqref{a3}.
\begin{theorem}\label{thm1}
 For $q\in (3, 6)$ and positive constant $\bar\rho$, assume that the initial data $(\rho_0, u_0, \theta_0\ge 0, H_0)$ satisfies
\begin{align}\label{qq1}
& 0\le \inf\rho_0(x)\le \sup\rho_0(x)\le \bar{\rho},\ \rho_0\in L^1\cap H^1\cap W^{1, q}, \\ \label{1.6}
& \nabla u_0\in H^1,\ \nabla\theta_0\in H^1, \ H_0\in H^2,
\end{align}
and the compatibility conditions
\begin{align}\label{qqw1}
& -\mu\Delta u_0-(\lambda+\mu)\nabla\divv u_0+\nabla p_0-\curl H_0\times H_0=\sqrt{\rho_0}g_1,\\
& \kappa\Delta\theta_0+\lambda(\divv u_0)^2+\frac{\mu}{2}|\nabla u_0+(\nabla u_0)^{tr}|^2+\eta|\curl H_0|^2=\sqrt{\rho_0}g_2,\label{qqw}
\end{align}
with $g_1, g_2\in L^2$.
There exists a small positive
constant $\varepsilon_0$ depending only on  $\mu$, $\lambda$, $c_v$, $\eta$, $R$, $\kappa$, $\bar{\rho}$, and the initial data such that if
\begin{align}\label{1.9}
\int\Big(\frac12\rho_0|u_0|^2+c_v\rho_0\theta_0+\frac12|H_0|^2\Big)dx
\triangleq\mathbb{E}_0\le \varepsilon_0,
\end{align}
then the problem \eqref{a1}--\eqref{a3} has a unique global strong solution $(\rho\geq0, u, \theta\geq0, H)$ satisfying
\begin{align}\label{trr}
\begin{cases}
\rho \in L^\infty\big([0,\infty);L^1\cap L^\infty\big)\cap C\big([0, \infty); L^1\cap H^1\cap W^{1, q}\big),\\
(u, \theta)\in C\big([0, \infty); D_0^1\cap D^2\big)\cap L_{loc}^2\big([0, \infty); D^{2, q}\big),\\
(u_t, \theta_t)\in L_{loc}^2\big([0, \infty); D_0^1\big), \ \big(\sqrt{\rho}u_t, \sqrt{\rho}\theta_t\big)\in L_{loc}^\infty\big([0, \infty); L^2\big),\\
H\in L^\infty\big([0,\infty);H^2\big)\cap C\big([0, \infty); H^2\big)\cap L^2\big([0, \infty); H^3\big), \\
H_t\in C\big([0, \infty); L^2\big)\cap L^2\big([0, \infty); H^1\big).
\end{cases}
\end{align}
\end{theorem}

\begin{remark}
It follows from H{\"o}lder's inequality, Sobolev's inequality, \eqref{qq1}, and \eqref{1.6} that
\begin{align*}
\int\big(\rho_0|u_0|^2+\rho_0\theta_0\big)dx
& \leq \|\rho_0\|_{L^\frac32}\|u_0\|_{L^6}^2
+\|\sqrt{\rho_0}\|_{L^2}\|\sqrt{\rho_0}\theta_0\|_{L^2} \\
& \leq C\|\rho_0\|_{L^\frac32}\|\nabla u_0\|_{L^2}^2
+\|\rho_0\|_{L^1}^{\frac12}\|\rho_0\|_{L^\frac32}^{\frac12}\|\theta_0\|_{L^6} \\
& \leq C\|\rho_0\|_{L^1}^{\frac23}\bar{\rho}^{\frac13}\|\nabla u_0\|_{L^2}^2
+C\|\rho_0\|_{L^1}^{\frac56}\bar{\rho}^{\frac16}\|\nabla \theta_0\|_{L^2}\\
& \leq C.
\end{align*}
Thus, $\mathbb{E}_0$ is finite.
\end{remark}

\begin{remark}
Compared with our previous work \cite{LZ20} where the global strong solution was established provided $\|\rho_0\|_{L^\infty}+\|H_0\|_{L^3}$ is suitably small, although our Theorem \ref{thm1} has small energy, yet whose oscillations could be arbitrarily large. Moreover, there is no additional restriction on the viscosity coefficients $\mu$ and $\lambda$ except the physical condition \eqref{a0}.
\end{remark}


As we know, the issue of global existence and uniqueness relies crucially on whether or not one can obtain global bounds on the solutions.
We make some techniques developed by the authors in \cite{L21,LZ20} to give the proof of Theorem \ref{thm1} and the key point in our analysis is to derive the time-independent upper bound of the density in Proposition \ref{p1}.
However, compared with the compressible non-isentropic Navier-Stokes system considered in \cite{L21}, due to the strong coupling and interplay interaction between the fluid motion and the magnetic field, the crucial techniques of proofs in \cite{L21} cannot be adapted directly to the situation treated
here. Moreover, unlike \cite{LZ20}, it is impossible to get the
smallness of $\|H\|_{L^3}$ which plays a critical role in dealing with the coupling terms $\curl(u\times H)$ and $\curl H\times H$.
To overcome these difficulties, some new ideas and observations are needed.

As pointed out in \cite{LZ20}, the basic energy estimate provides us
useless information on dissipation estimates of $u$, $\theta$, and $H$. To overcome this difficulty, by introducing a small time $t_0$ to be determined, we recover the crucial dissipation estimate of the form
$\int_0^T\big(\|\nabla u\|_{L^2}^2+\|\nabla \theta\|_{L^2}^2+\|\nabla H\|_{L^2}^2\big)dt$ in terms of a power function of $\mathbb{E}_0$ (see Lemma \ref{l31}) instead of $\|\rho_0\|_{L^\infty}+\|H_0\|_{L^3}$ used in \cite{LZ20}. Moreover, we should point out that the method used here to derive the $L^\infty(0,T;L^2)$ estimates of $\nabla u$ and $\nabla H$ is different from that of in \cite{LZ20}. In order to obtain this control, motivated by \cite{L21}, the key observation is to treat the temperature simultaneously by introducing a cut-off function $\sigma(t)\triangleq\min\{t,t_0\}$ (see Lemma \ref{l33}) instead of $E\triangleq\frac{|u|}{2}+c_v\theta$, the sum of specific kinetic energy and specific internal energy.
Next, combining these estimates altogether, we adopt some ideas
introduced by \cite{D1997,LZ20} to establish the time-independent upper bound of the density. Finally, having the time-independent and time-dependent
\textit{a priori} estimates at hand, together with the local well-posedness theory, we succeed in showing Theorem \ref{thm1} by the continuation method.

Our second main result is about algebraic decay properties of the solution showed in Theorem \ref{thm1}.
\begin{theorem}\label{thm2}
Let $(\rho, u, \theta, H)$ be the solution obtained by Theorem \ref{thm1}, then there exists another constant $\varepsilon_0'\in (0, \varepsilon_0]$
such that if
\begin{align*}
\mathbb{E}_0\le \varepsilon_0',
\end{align*}
then, for large time $t$, the following decay estimates hold true
\begin{align}\label{1.11}
&\|\sqrt{\rho}u\|_{L^2}+\|\sqrt{\rho}\theta\|_{L^2}
+\|\nabla u\|_{L^2}+\|\nabla\theta\|_{L^2}+\|\sqrt{\rho}\dot{u}\|_{L^2}\le Ct^{-\frac14},\\ \label{1.12}
&\|H_t\|_{L^2}+\|\nabla H\|_{H^1}
+\|\sqrt{\rho}\dot{\theta}\|_{L^2}+\|\nabla^2\theta\|_{L^2}\le Ct^{-\frac12},
\end{align}
where the constant $C$ depends only on $\mu$, $\lambda$, $\eta$, $c_v$, $\kappa$, $\bar{\rho}$, and the initial data.
\end{theorem}


\begin{remark}
Let $F$ and $w$ (see \eqref{f3}) be the effective viscous flux and the vorticity, respectively. Then it follows from \eqref{w5}, Sobolev's inequality, \eqref{trr}, \eqref{1.11}, and \eqref{1.12} that, for large time $t$,
\begin{align}\label{1.13}
\|\nabla F\|_{L^2}+\|\nabla w\|_{L^2}
& \le C\big(\|\rho\dot{u}\|_{L^2}+\||H||\nabla H|\|_{L^2}\big) \notag \\
& \leq C\|\rho\|_{L^\infty}^{\frac12}\|\sqrt{\rho}\dot{u}\|_{L^2}
+C\|H\|_{L^\infty}\|\nabla H\|_{L^2} \notag \\
& \leq C\|\rho\|_{L^\infty}^{\frac12}\|\sqrt{\rho}\dot{u}\|_{L^2}
+C\|H\|_{H^2}\|\nabla H\|_{L^2}
\notag \\
& \leq C\|\sqrt{\rho}\dot{u}\|_{L^2}+C\|\nabla H\|_{L^2} \notag \\
& \leq Ct^{-\frac14}.
\end{align}
This combined with \eqref{w7} and \eqref{1.11} yields that, for large time $t$,
\begin{align}\label{1.14}
\|\nabla u\|_{L^6}
\le C\big(\|\rho\dot{u}\|_{L^2}+\||H||\nabla H|\|_{L^2}
+\|\nabla\theta\|_{L^2}\big)\leq Ct^{-\frac14}.
\end{align}
Furthermore, we derive from H\"older's inequality, Sobolev's inequality, \eqref{1.14}, \eqref{1.11}, and \eqref{1.12} that, for any $q\in[2,6]$ and large time $t$,
\begin{align*}
& \|\nabla u\|_{L^q} \leq \|\nabla u\|_{L^2}^{\frac{6-q}{2q}}
\|\nabla u\|_{L^6}^{\frac{3q-6}{2q}}\leq Ct^{-\frac14},\\
& \|\nabla\theta\|_{L^q} \leq \|\nabla\theta\|_{L^2}^{\frac{6-q}{2q}}
\|\nabla\theta\|_{L^6}^{\frac{3q-6}{2q}}
\leq C\|\nabla\theta\|_{L^2}^{\frac{6-q}{2q}}
\|\nabla^2\theta\|_{L^2}^{\frac{3q-6}{2q}}\leq Ct^{-\frac{5q-6}{8q}}.
\end{align*}
\end{remark}

\begin{remark}
Although the large-time asymptotic decay rate of the gradient of the vorticity is shown in \eqref{1.13}, yet we should point out that whether the second
derivatives of the velocity field decay or not remains unknown.
The main difficulty lies in deriving time-independent bounds on the gradient of density (see \eqref{t91}), which in turn effects the estimate of the velocity (see \eqref{t87}).
\end{remark}

\begin{remark}
In \cite{L21} Liang established the exponential decay-in-time property of the solution of full compressible Navier-Stokes equations which is quite different from Theorem \ref{thm2} for the non-isentropic MHD equations. This means that the magnetic field acts as some significant roles on the large time behaviors of the velocity and the temperature.
\end{remark}

The proof of Theorem \ref{thm2} is based on delicate energy estimates. Firstly, due to the structure of \eqref{a1}, we find that we can only obtain \begin{align*}
\frac{d}{dt}\big(\|\sqrt{\rho}u\|_{L^2}^2+\|\sqrt{\rho}\theta\|_{L^2}^2\big)
+\|\nabla u\|_{L^2}^2+\|\nabla \theta\|_{L^2}^2
\le C\|\nabla H\|_{L^2}\|\nabla H\|_{H^1}^2
\end{align*}
for small initial energy (see \eqref{x4}), thus we need to show the time-decay rate of the gradient of the magnetic field in order to obtain the initial layer decay estimates of the velocity and the temperature. Fortunately, by analyzing \eqref{a1}$_4$, we derive that $\|\nabla H\|_{L^2}$ decays with the rate of $t^{-\frac12}$ (see \eqref{dd1}). Next, owing to very slowly time-decay estimates on the initial layer analysis,
we fail to establish the decay rates of $\|\nabla u\|_{L^2}$ and $\|\nabla \theta\|_{L^2}$ directly and should require the initial energy to be much more small (see \eqref{5.15}). Finally, the decay estimates on the second derivatives of the magnetic field and the temperature
can be obtained by applying $L^2$-estimates for the elliptic system and uniformly bounds we have just derived (see Step 6 in Section \ref{sec5}).

The rest of the paper is organized as follows. In Section 2, we recall some known facts and elementary inequalities which will be used later. Section 3 is devoted to the global \textit{a priori} estimates. The proof of Theorem \ref{thm1} will be done in Section 4. Finally, we give the proof of Theorem \ref{thm2} in Section 5.

\section{Preliminaries}
In this section, we recall some known facts and elementary inequalities which will be used frequently later.

We begin with the local existence of a unique strong solution with vacuum to the problem \eqref{a1}--\eqref{a3}, whose proof can be found in \cite{FY09}.
\begin{lemma}\label{l21}
Under the conditions \eqref{qq1}--\eqref{qqw}, there is a small $T_*>0$ such that the Cauchy problem \eqref{a1}--\eqref{a3} admits a unique solution $(\rho, u, \theta, H)$
over $\mathbb{R}^3\times(0, T_*]$ satisfying for some constant $M_*>1$ depending on $T_*$ and the initial data
\begin{align}\label{w1}
&\sup_{0\le t\le T_*}\big(\|\nabla u\|_{H^1}^2+\|\nabla\theta\|_{H^1}^2+\|H\|_{H^2}^2
+\|\sqrt{\rho}\dot{u}\|_{L^2}^2+\|\sqrt{\rho}\dot{\theta}\|_{L^2}^2+\|H_t\|_{L^2}^2\big)\nonumber\\
&+\int_0^{T_*}\big(\|\nabla\dot{u}\|_{L^2}^2+\|\nabla H_t\|_{L^2}^2+\|\nabla\dot{\theta}\|_{L^2}^2\big)dt\le M_*.
\end{align}
\end{lemma}

Next, the following well-known Gagliardo-Nirenberg inequality (see \cite[Chapter 10, Theorem 1.1]{D2016}) will be
used frequently later.
\begin{lemma}\label{l22}
Assume that $f\in D^{1, m}\cap L^r$ with $m,r\geq1$, then there exists a constant $C$ depending only on $q$, $m$, and $r$ such that
\begin{align}
\|f\|_{L^q}\le C\|\nabla f\|_{L^m}^\vartheta\|f\|_{L^r}^{1-\vartheta},
\end{align}
where $\vartheta=\big(\frac{1}{r}-\frac{1}{q}\big)/
\big(\frac{1}{r}-\frac{1}{m}+\frac13\big)$ and the admissible range of $q$ is the following:
\begin{itemize}
\item if $m<3$, then $q$ is between $r$ and $\frac{3m}{3-m}$;
\item if $m=3$, then $q\in [r, \infty)$;
\item if $m>3$, then $q\in [r, \infty]$.
\end{itemize}
\end{lemma}

Next,  the material derivatives of $f$, the effective viscous flux $F$,
and the vorticity $w$ are defined as follows
\begin{align}\label{f3}
\dot{f}\triangleq u_t+u\cdot\nabla f, \ \
F\triangleq(2\mu+\lambda)\divv u-p-\frac12|H|^2, \ \ w\triangleq\curl u.
\end{align}
Then one derives the following elliptic equations from $\eqref{a1}_2$
by taking the operators ${\rm div}$ and ${\rm curl}$, respectively.
\begin{align}\label{s22}
\Delta F=\divv(\rho\dot{u})-\divv\divv(H\otimes H),\ \
\mu\Delta w=\curl \big(\rho\dot{u}+\divv(H\otimes H)\big).
\end{align}
We now state some elementary $L^q$-estimates of the elliptic system \eqref{s22}, whose proof can be done by similar methods as those in \cite{LXZ13,L21}.
\begin{lemma}
Let $(\rho, u, \theta, H)$ be a smooth solution of \eqref{a1}--\eqref{a3}. Then for any $q\in[2,6]$, there exists a positive constant $C$ depending only on $\mu$, $\lambda$, $\eta$, and $q$ such that
\begin{align}
&\|\nabla F\|_{L^q}+\|\nabla w\|_{L^q}\le C\big(\|\rho\dot{u}\|_{L^q}+\||H||\nabla H|\|_{L^q}\big),\label{w5}\\
&\|\nabla u\|_{L^6}\le C\big(\|\rho\dot{u}\|_{L^2}+\|\nabla\theta\|_{L^2}
+\||H||\nabla H|\|_{L^2}\big).\label{w7}
\end{align}
\end{lemma}

Finally, the following Beale-Kato-Majda type inequality (see \cite[Lemma 2.3]{HLX20112}) will be used to estimate $\|\nabla u\|_{L^\infty}$.
\begin{lemma}\label{l24}
Suppose that $q\in (3, \infty)$, then there exists a constant $C$ depending only on $q$ such that for any $\nabla v\in L^2\cap W^{1, q}$,
\begin{align}
\|\nabla v\|_{L^\infty}\le C\big(\|\divv v\|_{L^\infty}+\|\curl v\|_{L^\infty}\big)\ln\big(e+\|\nabla^2v\|_{L^q}\big)+C\big(1+\|\nabla v\|_{L^2}\big).
\end{align}
\end{lemma}

\section{\textit{A priori} estimates}
In this section, we will establish some necessary {\it a priori} bounds for
smooth solutions to the Cauchy problem \eqref{a1}--\eqref{a3} to extend the local strong
solutions guaranteed by Lemma \ref{l21}. Thus, let $T>0$ be a fixed time and
$(\rho, u, \theta, H)$ be the smooth solution to \eqref{a1}--\eqref{a3} in $\mathbb{R}^3\times(0, T]$ with initial data
$(\rho_0, u_0, \theta_0, H_0)$ satisfying \eqref{qq1}--\eqref{qqw}. In what follows, $C, C_i$, and $c_i\ (i=1, 2, \cdots)$ denote generic positive
constants which rely only on $R$, $c_v$, $\mu$, $\lambda$, $\kappa$, $\eta$, $\bar{\rho}$, $M_*$, and the initial values, and $C(\alpha)$
is used to emphasize the dependence of $C$ on $\alpha$.

First of all, if we multiply $\eqref{a1}_1$ by a cut-off function and then use a standard limit procedure, we can obtain that (see also \cite[Lemma 3.1]{WZ17})
\begin{align}\label{lz3.1}
\|\rho(t)\|_{L^1}=\|\rho_0\|_{L^1}, \quad 0\leq t\leq T.
\end{align}
Moreover, applying standard maximum principle (see \cite[p. 43]{F2004}) to \eqref{a1} along with $\rho_0,\theta_0\geq0$ shows
\begin{equation}\label{lz3.2}
\inf_{\mathbb{R}^3\times[0,T]}\rho(x,t)\geq0,\ \ \inf_{\mathbb{R}^3\times[0,T]}\theta(x,t)\geq0.
\end{equation}
Meanwhile, the basic energy estimate gives that
\begin{equation}\label{lz3.3}
\int\Big(\frac12\rho|u|^2+c_v\rho\theta+\frac12|H|^2\Big)dx
=\int\Big(\frac12\rho_0|u_0|^2+c_v\rho_0\theta_0+\frac12|H_0|^2\Big)dx,
\quad 0\leq t\leq T.
\end{equation}

\subsection{Time-independent estimates}
Denote
\begin{align*}
A_1(T)\triangleq &\sup_{0\le t\le T}\big(\|\nabla u\|_{L^2}^2+\|\nabla H\|_{L^2}^2+\|\nabla\theta\|_{L^2}^2\big),\\
A_2(T)\triangleq &\int_0^T\big(\|\sqrt{\rho}\dot{u}\|_{L^2}^2+\|H_t\|_{L^2}^2
+\|\nabla^2 H\|_{L^2}^2+\|\nabla\theta\|_{H^1}^2\big)dt.
\end{align*}

\begin{proposition}\label{p1}
There exists some positive constant $\varepsilon_0$ depending only on
$\mu$, $\lambda$, $c_v$, $\eta$, $R$, $\kappa$, $\bar{\rho}$, and the initial data such that if $(\rho, u, \theta, H)$ is the smooth solution of \eqref{a1}--\eqref{a3} on $\mathbb{R}^3\times(0, T]$ satisfying
\begin{align}\label{ee1}
\sup_{0\le t\le T}\|\rho\|_{L^\infty}\le 2\bar{\rho}, \ A_1(T)\le 2 M_*, \ A_2(T)\le 2\mathbb{E}_0^\frac{2}{3},
\end{align}
then the following estimates hold
\begin{align}\label{ee2}
\sup_{0\le t\le T}\|\rho\|_{L^\infty}\le \frac{3}{2}\bar{\rho}, \ A_1(T)\le \frac{3}{2}M_*, \ A_2(T)\le \mathbb{E}_0^\frac{2}{3},
\end{align}
provided that
$\mathbb{E}_0\le \varepsilon_0.$
\end{proposition}

Before proving Proposition \ref{p1}, we establish some necessary \textit{a priori} estimates, see Lemmas \ref{l31}--\ref{l34} below.
\begin{lemma}\label{l31}
Let $(\rho, u, \theta, H)$ be the smooth solution of \eqref{a1}--\eqref{a3} satisfying \eqref{ee1}, then we have
\begin{align}\label{f5}
&\sup_{0\le t\le T}\big(\|\sqrt{\rho}u\|_{L^2}^2+\|\sqrt{\rho}\theta\|_{L^2}^2
+\|H\|_{L^2}^2\big)+\int_0^T\big(\|\nabla u\|_{L^2}^2+\|\nabla\theta\|_{L^2}^2
+\|\nabla H\|_{L^2}^2\big)dt\nonumber\\
&\le C\Big(t_0+\mathbb{E}_0^\frac{4}{5}\Big)+C\mathbb{E}_0^\frac{2}{3}\sup_{t_0\le t\le T}\big(\|\nabla u\|_{L^2}^2+\|\nabla\theta\|_{L^2}^2
+\|\nabla H\|_{L^2}^2\big),
\end{align}
where $t_0\in (0, T_*)$ is to be determined and $T_*$ is taken from Lemma \ref{l21}.
\end{lemma}
{\it Proof.}
1. Due to ${\rm div}H=0$, direct calculation gives that
\begin{align*}
&\curl(u\times H)=H\cdot\nabla u-u\cdot\nabla H-H\divv u,\\
&\curl H\times H=H\cdot\nabla H-\frac12\nabla|H|^2.
\end{align*}
Multiplying $\eqref{a1}_2$ by $u$ and $\eqref{a1}_4$ by $H$, respectively, and integrating the resulting equations over $\mathbb{R}^3$ and adding them together, we obtain from H\"older's inequality, Sobolev's inequality, \eqref{lz3.1}, and \eqref{ee1} that
\begin{align*}
&\frac{1}{2}\frac{d}{dt}\int\big(\rho|u|^2+|H|^2\big)dx
+\int\Big[\mu|\nabla u|^2+(\mu+\lambda)(\divv u)^2
+\eta|\nabla H|^2\Big]dx\nonumber\\
&=\int R\rho\theta\divv udx\nonumber\\
&\le (\mu+\lambda)\|\divv u\|_{L^2}^2
+C\|\rho\|_{L^\infty}\|\sqrt{\rho}\theta\|_{L^2}^2\nonumber\\
&\le (\mu+\lambda)\|\divv u\|_{L^2}^2
+C(\bar{\rho})\|\rho\|_{L^\frac{3}{2}}\|\theta\|_{L^6}^2\nonumber\\
&\le (\mu+\lambda)\|\divv u\|_{L^2}^2
+C\|\rho\|_{L^1}^\frac23\|\rho\|_{L^\infty}^\frac{1}{3}
\|\nabla\theta\|_{L^2}^2\nonumber\\
&\le (\mu+\lambda)\|\divv u\|_{L^2}^2+c_1\|\nabla\theta\|_{L^2}^2,
\end{align*}
for some constant $c_1$ depending only $\mu,\lambda,\bar\rho$, and $\|\rho_0\|_{L^1}$. Thus, we get that
\begin{align}\label{f6}
&\frac{d}{dt}\big(\|\sqrt{\rho}u\|_{L^2}^2+\|H\|_{L^2}^2\big)
+2\mu\|\nabla u\|_{L^2}^2+2\eta\|\nabla H\|_{L^2}^2
\leq 2c_1\|\nabla\theta\|_{L^2}^2.
\end{align}

2. Multiplying $\eqref{a1}_3$ by $\theta$ and integration by parts, we obtain from Sobolev's inequality that
\begin{align*}
&\frac{c_v}{2}\frac{d}{dt}\int\rho\theta^2dx+\kappa\int|\nabla\theta|^2dx
\nonumber\\
&=\int\theta\Big(\lambda(\divv u)^2+\frac{\mu}{2}|\nabla u+(\nabla u)^{tr}|^2-R\rho\theta\divv u\Big)dx+\eta\int\theta\curl H\cdot\curl H dx\nonumber\\
&\le \frac{\mu\kappa}{6c_1}\int|\nabla u|^2dx+C\int\theta^2\big(|\nabla u|^2+\rho\theta^2\big)dx+C\int\theta|\nabla^2H|H|dx+C\int|\nabla\theta||\nabla H|H|dx\nonumber\\
&\le \frac{\mu\kappa}{6c_1}\|\nabla u\|_{L^2}^2+C\|\theta\|_{L^\infty}^2\big(\|\nabla u\|_{L^2}^2+\|\sqrt{\rho}\theta\|_{L^2}^2\big)
+C\big(\|\theta\|_{L^6}\|\nabla^2 H\|_{L^2}+\|\nabla\theta\|_{L^2}\|\nabla H\|_{L^6}\big)\|H\|_{L^3}
\nonumber\\
&\le \frac{\mu\kappa}{6c_1}\|\nabla u\|_{L^2}^2+C\|\theta\|_{L^\infty}^2\big(\|\nabla u\|_{L^2}^2+\|\sqrt{\rho}\theta\|_{L^2}^2\big)
+C\|\nabla\theta\|_{L^2}\|\nabla^2H\|_{L^2}\|H\|_{L^2}^\frac12\|\nabla H\|_{L^2}^\frac12\nonumber\\
&\le \frac{\mu\kappa}{6c_1}\|\nabla u\|_{L^2}^2+\frac{\kappa}{2}\|\nabla\theta\|_{L^2}^2+C\|\theta\|_{L^\infty}^2
\big(\|\nabla u\|_{L^2}^2+\|\nabla\theta\|_{L^2}^2\big)
+C\|\nabla^2H\|_{L^2}^2\big(\|H\|_{L^2}^2+\|\nabla H\|_{L^2}^2\big),
\end{align*}
which implies that
\begin{align}\label{f7}
&c_v\frac{d}{dt}\|\sqrt{\rho}\theta\|_{L^2}^2
+\kappa\|\nabla\theta\|_{L^2}^2\nonumber\\
&\le \frac{\mu\kappa}{3c_1}\|\nabla u\|_{L^2}^2+C\|\theta\|_{L^\infty}^2\big(\|\nabla u\|_{L^2}^2+\|\nabla\theta\|_{L^2}^2\big)+C\|\nabla^2H\|_{L^2}^2
\big(\|H\|_{L^2}^2+\|\nabla H\|_{L^2}^2\big).
\end{align}
Adding \eqref{f7} multiplied by $\frac{3c_1}{\kappa}$ to \eqref{f6}, we deduce that
\begin{align}\label{f9}
&\frac{d}{dt}\big(\|\sqrt{\rho}u\|_{L^2}^2+\|H\|_{L^2}^2
+\|\sqrt{\rho}\theta\|_{L^2}^2\big)+\|\nabla u\|_{L^2}^2+\|\nabla\theta\|_{L^2}^2+\|\nabla H\|_{L^2}^2\nonumber\\
&\le C\|\theta\|_{L^\infty}^2\big(\|\nabla u\|_{L^2}^2+\|\nabla\theta\|_{L^2}^2\big)
+C\|\nabla^2H\|_{L^2}^2\big(\|H\|_{L^2}^2+\|\nabla H\|_{L^2}^2\big).
\end{align}

3. Applying the Gagliardo-Nirenberg inequality
\begin{align}\label{zx}
\|\theta\|_{L^\infty}\le C\|\nabla \theta\|_{L^2}^\frac12\|\nabla^2 \theta\|_{L^2}^\frac12,
\end{align}
we derive from \eqref{ee1} that
\begin{align}\label{g12}
\int_0^T\|\theta\|_{L^\infty}^2dt\le C\int_0^T\big(\|\nabla\theta\|_{L^2}^2+\|\nabla^2\theta\|_{L^2}^2\big)dt\le C\mathbb{E}_0^\frac23.
\end{align}
By H{\"o}lder's inequality, \eqref{1.6}, Sobolev's inequality (see \cite[Theorem]{T1976}), and \eqref{1.9}, we have
\begin{align}\label{3.11}
\|\sqrt{\rho_0}\theta_0\|_{L^2}^2
& \le \|\rho_0\theta_0^\frac45\|_{L^{\frac54}}\|\theta_0^\frac65\|_{L^5} \notag \\
& \le \|\rho_0\|_{L^\infty}^{\frac15}\|\rho_0\theta_0\|_{L^1}^\frac45
\|\theta_0\|_{L^6}^\frac65 \notag \\
& \le \bar{\rho}^{\frac15}c_v^{-\frac45}\|c_v\rho_0\theta_0\|_{L^1}^\frac45
\frac{1}{3^{\frac35}}\Big(\frac2\pi\Big)^{\frac45}
\|\nabla\theta_0\|_{L^2}^\frac65 \notag \\
& \le \bar{\rho}^{\frac15}\frac{1}{3^{\frac35}}\Big(\frac{2}{c_v\pi}\Big)^{\frac45}
\|\nabla\theta_0\|_{L^2}^\frac65\mathbb{E}_0^\frac45.
\end{align}
Integrating \eqref{f9} with respect to $t$ over $[0,T]$, we derive from \eqref{1.9}, \eqref{g12}, and \eqref{3.11} that
\begin{align}\label{zx2}
&\sup_{0\le t\le T}\big(\|\sqrt{\rho}u\|_{L^2}^2+\|\sqrt{\rho}\theta\|_{L^2}^2
+\|H\|_{L^2}^2\big)+\int_0^T\big(\|\nabla u\|_{L^2}^2+\|\nabla\theta\|_{L^2}^2
+\|\nabla H\|_{L^2}^2\big)dt\nonumber\\
&\le 2\mathbb{E}_0+\bar{\rho}^{\frac15}\frac{1}{3^{\frac35}}\Big(\frac{2}{c_v\pi}\Big)^{\frac45}
\|\nabla\theta_0\|_{L^2}^\frac65\mathbb{E}_0^\frac{4}{5} \notag \\
& \quad +C\sup_{0\le t\le t_0}\big(\|\nabla u\|_{L^2}^2+\|\nabla\theta\|_{L^2}^2
+\|\nabla H\|_{L^2}^2\big)
\int_0^{t_0}\big(\|\theta\|_{L^\infty}^2+\|\nabla^2H\|_{L^2}^2\big)dt\nonumber\\
&\quad+C\sup_{t_0\le t\le T}\big(\|\nabla u\|_{L^2}^2+\|\nabla\theta\|_{L^2}^2
+\|\nabla H\|_{L^2}^2\big)\int_{t_0}^T\big(\|\theta\|_{L^\infty}^2+\|\nabla^2H\|_{L^2}^2\big)dt\nonumber\\
&\quad+C\sup_{0\le t\le T}\|H\|_{L^2}^2\int_0^T\|\nabla^2 H\|_{L^2}^2dt\nonumber\\
&\le C\Big(t_0+\mathbb{E}_0^\frac{4}{5}\Big)+C\mathbb{E}_0^\frac{2}{3}
\sup_{t_0\le t\le T}\big(\|\nabla u\|_{L^2}^2+\|\nabla\theta\|_{L^2}^2
+\|\nabla H\|_{L^2}^2\big)+C_1\mathbb{E}_0^\frac{2}{3}\sup_{0\le t\le T}\|H\|_{L^2}^2,
\end{align}
provided that
\begin{align*}
\mathbb{E}_0\leq \frac{\bar\rho}{2}\frac{1}{3^3}\Big(\frac{1}{c_v\pi}\Big)^4
\|\nabla\theta_0\|_{L^2}^6
=\frac{\bar\rho\|\nabla\theta_0\|_{L^2}^6}{54c_v^4\pi^4},
\end{align*}
where in the last inequality of \eqref{zx2} we have used
\begin{align*}
& C\sup_{0\le t\le t_0}\big(\|\nabla u\|_{L^2}^2+\|\nabla\theta\|_{L^2}^2
+\|\nabla H\|_{L^2}^2\big)
\int_0^{t_0}\big(\|\theta\|_{L^\infty}^2+\|\nabla^2H\|_{L^2}^2\big)dt \\
&\leq Ct_0\sup_{0\le t\le t_0}\big(\|\nabla u\|_{L^2}^2+\|\nabla\theta\|_{L^2}^2
+\|\nabla H\|_{L^2}^2\big)
\sup_{0\le t\le t_0}\big(\|\theta\|_{L^\infty}^2+\|\nabla^2H\|_{L^2}^2\big) \\
&\leq Ct_0,
\end{align*}
due to \eqref{w1} and \eqref{zx}.
Thus, \eqref{f5} follows from \eqref{zx2}
provided that
$$\mathbb{E}_0\le \varepsilon_1\triangleq\min
\left\{\frac{\bar\rho\|\nabla\theta_0\|_{L^2}^6}{54c_v^4\pi^4}, \frac{1}{(2C_1)^\frac32}\right\}.$$
This finishes the proof of Lemma \ref{l31}.
\hfill $\Box$

\begin{lemma}\label{l32}
Let $(\rho, u, \theta, H)$ be the smooth solution of \eqref{a1}--\eqref{a3} satisfying \eqref{ee1}, then it holds that
\begin{align}\label{f12}
&\frac{d}{dt}\int\Big[\frac{\mu}{2}|\nabla u|^2+\frac{\kappa\delta}{2}|\nabla\theta|^2+\eta|\nabla H|^2
+\Big(\frac{2c_2}{\mu}+1\Big)\delta^\frac32\rho|\dot{u}|^2
+2\delta^\frac32|H_t|^2\Big]dx\nonumber\\
&\quad+\int\Big(\rho|\dot{u}|^2+\frac12|H_t|^2+\frac{\eta^2}{4}|\nabla^2H|^2
+\frac{c_v\delta}{4}\rho|\dot{\theta}|^2
+\mu \delta^\frac32|\nabla \dot{u}|^2+\frac{\eta\delta^\frac32}{2}|\nabla H_t|^2\Big)dx\nonumber\\
&\le \frac{d}{dt}\int\Big(R\rho\theta\divv u+\frac12|H|^2\divv u-H\cdot\nabla u\cdot H
+\lambda\theta\delta(\divv u)^2+\frac{\mu}{2}\theta\delta|\nabla u+(\nabla u)^{tr}|^2\Big)dx\nonumber\\
&\quad+\frac{d}{dt}\int\delta\theta|\curl H|^2dx
+\delta^{-\frac32}(\|\theta\nabla u\|_{L^2}^2+\|\nabla u\|_{L^2}^2+\|\nabla u\|_{L^2}^4\|\nabla\theta\|_{L^2}^2)
\nonumber\\
&\quad+\delta^\frac32C\|\nabla u\|_{L^4}^4+\delta^\frac52\|\nabla^2H\|_{L^2}^4
+C\|\nabla H\|_{L^2}^4+C\|\nabla u\|_{L^2}^4\|\nabla H\|_{L^2}^2
+C\|H\|_{L^2}\|\nabla H\|_{L^2}^3\nonumber\\[3pt]
&\quad+C\delta\|\theta\|_{L^\infty}^2\|\nabla H\|_{L^2}^2,
\end{align}
where $\delta$ is a small positive constant and to be determined.
\end{lemma}
{\it Proof.}
1. Multiplying $\eqref{a1}_2$ by $\dot{u}$ and then integrating the resulting equality over $\mathbb{R}^3$
lead to
\begin{align}\label{f13}
\int\rho|\dot{u}|^2dx&=-\int\dot{u}\cdot\nabla pdx+\mu\int\Delta u\cdot\dot{u}dx+(\mu+\lambda)\int\nabla\divv u\cdot\dot{u}dx\nonumber\\
&\quad+\int H\cdot\nabla H\cdot\dot{u}dx-\frac12\int\nabla|H|^2\cdot\dot{u}dx\triangleq\sum_{i=1}^5J_i.
\end{align}
By using the fact $\rho\dot{\theta}=(\rho\theta)_t+{\rm div}(\rho u\theta)$ due to \eqref{a1}$_1$, we infer from integration by parts and \eqref{ee1} that, for any given $\delta>0$,
\begin{align}\label{f14}
J_1&=\int p\divv u_tdx-\int(u\cdot\nabla u)\cdot\nabla pdx\nonumber\\
&=\frac{d}{dt}\int R\rho\theta\divv udx-R\int\big[(\rho\theta)_t+{\rm div}(\rho\theta u)\big]\divv udx+\int R\rho\theta\partial_ju^k\partial_ku^jdx\nonumber\\
&\le \frac{d}{dt}\int R\rho\theta\divv udx
+\delta^\frac32\big(\|\sqrt{\rho}\dot{\theta}\|_{L^2}^2+\|\theta\nabla u\|_{L^2}^2\big)
+C\delta^{-\frac32}\|\nabla u\|_{L^2}^2.
\end{align}
Here and in what follows, we use the Einstein convention that the repeated indices denote the summation.
Integration by parts implies that
\begin{align}
J_2&=-\frac{\mu}{2}\frac{d}{dt}\|\nabla u\|_{L^2}^2-\mu\int\partial_iu^j\partial_i\big(u^k\partial_ku^j\big)dx\nonumber\\
&\le -\frac{\mu}{2}\frac{d}{dt}\|\nabla u\|_{L^2}^2-\mu\int\Big(\partial_iu^j\partial_iu^k
\partial_ku^j-\frac{1}{2}|\nabla u|^2\divv u\Big)dx\nonumber\\
&\le -\frac{\mu}{2}\frac{d}{dt}\|\nabla u\|_{L^2}^2+C\|\nabla u\|_{L^3}^3.
\end{align}
Similarly, we have
\begin{align}
J_3&\le -\frac{\mu+\lambda}{2}\frac{d}{dt}\|\divv u\|_{L^2}^2+C\|\nabla u\|_{L^3}^3.
\end{align}
Moreover, we deduce from $\eqref{a1}_4$, $\eqref{a1}_5$, integration by parts, H\"older's inequality, and Gagliardo-Nirenberg inequality that
\begin{align}\label{3.18}
J_4&=-\frac{d}{dt}\int H\cdot\nabla u\cdot Hdx+\int H_t\cdot\nabla u\cdot Hdx+\int H\cdot\nabla u\cdot H_tdx+\int H\cdot\nabla H\cdot(u\cdot\nabla u)dx\nonumber\\
&=-\frac{d}{dt}\int H\cdot\nabla u\cdot Hdx+\int(\Delta H-u\cdot\nabla H+H\cdot\nabla u-H\divv u)\cdot\nabla u\cdot Hdx\nonumber\\
&\quad+\int H\cdot\nabla u\cdot(\Delta H-u\cdot\nabla H+H\cdot\nabla u-H\divv u)dx-\int H^i\partial_iu^j\partial_ju^kH^kdx\nonumber\\
&\quad-\int H^iu^j\partial_i\partial_ju^kH^kdx\nonumber\\
&=-\frac{d}{dt}\int H\cdot\nabla u\cdot Hdx+\int(\Delta H+H\cdot\nabla u-H\divv u)\cdot\nabla u\cdot Hdx\nonumber\\
&\quad+\int H\cdot\nabla u\cdot(\Delta H+H\cdot\nabla u-H\divv u)dx-\int H^i\partial_iu^j\partial_ju^kH^kdx\nonumber\\
&\le -\frac{d}{dt}\int H\cdot\nabla u\cdot Hdx+\frac{\eta^2}{16}\|\nabla^2 H\|_{L^2}^2+C\||H||\nabla u|\|_{L^2}^2\nonumber\\
&\le  -\frac{d}{dt}\int H\cdot\nabla u\cdot Hdx+\frac{\eta^2}{16}\|\nabla^2 H\|_{L^2}^2
+C\delta^\frac32\|\nabla u\|_{L^4}^4+C\|H\|_{L^4}^4\nonumber\\
&\le  -\frac{d}{dt}\int H\cdot\nabla u\cdot Hdx+\frac{\eta^2}{16}\|\nabla^2 H\|_{L^2}^2
+C\delta^\frac32\|\nabla u\|_{L^4}^4+C\|H\|_{L^2}\|\nabla H\|_{L^2}^3.
\end{align}
Similarly to \eqref{3.18}, we get that
\begin{align}\label{f18}
J_5\le  \frac12\frac{d}{dt}\int|H|^2\divv udx+\frac{\eta^2}{16}\|\nabla^2 H\|_{L^2}^2
+C\delta^\frac32\|\nabla u\|_{L^4}^4+C\|H\|_{L^2}\|\nabla H\|_{L^2}^3.
\end{align}
Inserting \eqref{f14}--\eqref{f18} into \eqref{f13} leads to
\begin{align}\label{f19}
&\frac12\frac{d}{dt}\int\big[\mu|\nabla u|^2+(\mu+\lambda)
(\divv u)^2\big]dx+\int\rho|\dot{u}|^2dx\nonumber\\
&\le \frac{d}{dt}\int\Big(R\rho\theta\divv u+\frac12|H|^2\divv u-H\cdot\nabla u\cdot H
\Big)dx+\frac{\eta^2}{8}\|\nabla^2H\|_{L^2}^2\nonumber\\
&\quad+C\delta^\frac32\big(\|\sqrt{\rho}\dot{\theta}\|_{L^2}^2+\|\theta\nabla u\|_{L^2}^2+\|\nabla u\|_{L^4}^4\big)
+C\delta^{-\frac32}\|\nabla u\|_{L^2}^2+C\|H\|_{L^2}\|\nabla H\|_{L^2}^3.
\end{align}

2. It follows from $\eqref{a1}_4$, H\"older's inequality, Gagliardo-Nirenberg inequality, and Sobolev's inequality that
\begin{align}\label{f20}
\eta\frac{d}{dt}\|\nabla H\|_{L^2}^2+\|H_t\|_{L^2}^2+\eta^2\|\nabla^2 H\|_{L^2}^2
&=\int|H_t-\eta\Delta H|^2dx\nonumber\\
&=\int|H\cdot\nabla u-u\cdot\nabla H-H\divv u|^2dx\nonumber\\
&\le C\||\nabla u||H|\|_{L^2}^2+C\||u||\nabla H|\|_{L^2}^2\nonumber\\
&\le C\|H\|_{L^\infty}^2\|\nabla u\|_{L^2}^2+C\|u\|_{L^6}^2\|\nabla H\|_{L^3}^2\nonumber\\
&\le C\|\nabla u\|_{L^2}^2\|\nabla H\|_{L^2}\|\nabla^2 H\|_{L^2}\nonumber\\
&\le \frac{\eta^2}{8}\|\nabla^2H\|_{L^2}^2+C\|\nabla u\|_{L^2}^4\|\nabla H\|_{L^2}^2.
\end{align}

3. Multiplying $\eqref{a1}_3$ by $\dot{\theta}$ and integrating the resulting equation by parts yields that
\begin{align}\label{f21}
&\frac{\kappa}{2}\frac{d}{dt}\int|\nabla\theta|^2dx+c_v\int\rho|\dot{\theta}|^2dx\nonumber\\
&=\kappa\int u\cdot\nabla\theta\Delta\theta dx
-\int\rho\theta\divv u\dot{\theta}dx
+\lambda\int(\divv u)^2\dot{\theta}dx\nonumber\\
&\quad+\frac{\mu}{2}\int|\nabla u+(\nabla u)^{tr}|^2\dot{\theta}dx
+\eta\int|\curl H|^2\dot{\theta}dx\triangleq\sum_{i=1}^5N_i.
\end{align}
The standard $L^2$-estimate of elliptic equations to \eqref{a1}$_3$ together with \eqref{ee1} gives that
\begin{align}\label{f22}
\|\nabla^2\theta\|_{L^2}&\le C(\|\rho\dot{\theta}\|_{L^2}
+\|\rho\theta\divv u\|_{L^2}+\|\nabla u\|_{L^4}^2
+\|\nabla H\|_{L^4}^2)\nonumber\\
&\le C(\bar\rho)\|\sqrt{\rho}\dot{\theta}\|_{L^2}
+C(\bar\rho)\|\theta\nabla u\|_{L^2}+C\|\nabla u\|_{L^4}^2
+C\|\nabla H\|_{L^4}^2,
\end{align}
which combined with Sobolev's inequality, Gagliardo-Nirenberg inequality, and Young's inequality implies that
\begin{align}\label{z3.23}
N_1&\leq C\|u\|_{L^6}\|\nabla\theta\|_{L^3}\|\nabla^2\theta\|_{L^2}\nonumber\\
&\le C\|\nabla u\|_{L^2}\|\nabla\theta\|_{L^2}^\frac12\|\nabla^2\theta\|_{L^2}^\frac32\nonumber\\
&\le C\|\nabla u\|_{L^2}\|\nabla\theta\|_{L^2}^\frac12
\big(\|\sqrt{\rho}\dot{\theta}\|_{L^2}+\|\theta\nabla u\|_{L^2}+\|\nabla u\|_{L^4}^2+\|\nabla H\|_{L^4}^2\big)^\frac32\nonumber\\
&\le \delta^\frac12\big(\|\sqrt{\rho}\dot{\theta}\|_{L^2}^2+\|\theta\nabla u\|_{L^2}^2+\|\nabla u\|_{L^4}^4
+\|\nabla H\|_{L^4}^4\big)+C\delta^{-\frac32}\|\nabla u\|_{L^2}^4\|\nabla\theta\|_{L^2}^2\nonumber\\
&\le \delta^\frac12\big(\|\sqrt{\rho}\dot{\theta}\|_{L^2}^2
+\|\theta\nabla u\|_{L^2}^2+\|\nabla u\|_{L^4}^4\big)+
C\delta^{-\frac32}\|\nabla u\|_{L^2}^4\|\nabla\theta\|_{L^2}^2+\delta^\frac12\|\nabla H\|_{L^2}\|\nabla^2H\|_{L^2}^3\nonumber\\
&\le \delta^\frac12\big(\|\sqrt{\rho}\dot{\theta}\|_{L^2}^2+\|\theta\nabla u\|_{L^2}^2+\|\nabla u\|_{L^4}^4\big)+
C\delta^{-\frac32}\|\nabla u\|_{L^2}^4\|\nabla\theta\|_{L^2}^2
\notag \\
& \quad +\delta^\frac32\|\nabla^2H\|_{L^2}^4
+C\delta^{-\frac52}\|\nabla H\|_{L^2}^4.
\end{align}
By \eqref{ee1} and Cauchy-Schwarz inequality, we have
\begin{align}
N_2\le \delta^\frac{1}{2}\|\sqrt{\rho}\dot{\theta}\|_{L^2}^2+C\delta^{-\frac{1}{2}}\|\theta\nabla u\|_{L^2}^2.
\end{align}
Integration by parts together with \eqref{f3} leads to
\begin{align}\label{f25}
N_3&=\lambda\int(\divv u)^2\theta_tdx+\lambda\int(\divv u)^2(u\cdot\nabla\theta)dx\nonumber\\
&=\lambda\frac{d}{dt}\int(\divv u)^2\theta dx-2\lambda\int\theta\divv u\divv\dot{u}dx
+2\lambda\int\theta\divv u\divv(u\cdot\nabla u)dx+\lambda\int(\divv u)^2(u\cdot\nabla\theta)dx\nonumber\\
&=\lambda\frac{d}{dt}\int(\divv u)^2\theta dx-2\lambda\int\theta\divv u\divv\dot{u}dx
+2\lambda\int\theta\divv u\partial_iu^j\partial_ju^idx+\lambda\int u\cdot\nabla(\theta(\divv u)^2)dx\nonumber\\
&=\lambda\frac{d}{dt}\int(\divv u)^2\theta dx-2\lambda\int\theta\divv u\divv\dot{u}dx
+2\lambda\int\theta\divv u\partial_iu^j\partial_ju^idx-\lambda\int\theta(\divv u)^3dx\nonumber\\
&\le \lambda\frac{d}{dt}\int(\divv u)^2\theta dx
+\delta(\|\nabla u\|_{L^4}^4+\|\nabla\dot{u}\|_{L^2}^2)+C\delta^{-1}\|\theta\nabla u\|_{L^2}^2.
\end{align}
Similarly to \eqref{f25}, one gets
\begin{align}
N_4&\le \frac{\mu}{2}\frac{d}{dt}\int\theta|\nabla u+(\nabla u)^{tr}|^2
+\delta(\|\nabla u\|_{L^4}^4+\|\nabla\dot{u}\|_{L^2}^2)+C\delta^{-1}\|\theta\nabla u\|_{L^2}^2.
\end{align}
For the term $N_5$, we deduce from Gagliardo-Nirenberg, H\"older's, and Young's inequalities that
\begin{align}\label{z3.27}
N_5&=\int|\curl H|^2\theta_tdx+\int|\curl H|^2(u\cdot\nabla\theta)dx\nonumber\\
&=\frac{d}{dt}\int|\curl H|^2\theta dx-\int(|\curl H|^2)_t\theta dx
+\int|\curl H|^2(u\cdot\nabla\theta)dx\nonumber\\
&\le \frac{d}{dt}\int|\curl H|^2\theta dx
+C\|\theta\|_{L^\infty}\|\nabla H\|_{L^2}\|\nabla H_t\|_{L^2}+C\|\nabla H\|_{L^6}^2\|u\|_{L^6}\|\nabla\theta\|_{L^2}\nonumber\\
&\le \frac{d}{dt}\int|\curl H|^2\theta dx
+C\|\theta\|_{L^\infty}^2\|\nabla H\|_{L^2}^2+\delta^\frac12\|\nabla H_t\|_{L^2}^2+C\|\nabla u\|_{L^2}\|\nabla\theta\|_{L^2}\|\nabla^2H\|_{L^2}^2.
\end{align}
Inserting \eqref{z3.23}--\eqref{z3.27} into \eqref{f21}, we find that
\begin{align}\label{f28}
&\frac{\kappa}{2}\frac{d}{dt}\int|\nabla\theta|^2dx
+\big(c_v-2\delta^\frac12\big)\int\rho|\dot{\theta}|^2dx\nonumber\\
&\le \frac{d}{dt}\int\theta\Big(\lambda(\divv u)^2+\frac{\mu}{2}|\nabla u+(\nabla u)^{tr}|^2+
|\nabla H|^2\Big)dx\nonumber\\
&\quad+\delta^{-\frac32}\big(\|\theta\nabla u\|_{L^2}^2+\|\nabla u\|_{L^2}^4\|\nabla\theta\|_{L^2}^2\big)
+3\delta^\frac12\|\nabla u\|_{L^4}^4+2\delta\|\nabla\dot{u}\|_{L^2}^2\nonumber\\[3pt]
&\quad+C_2\|\nabla^2H\|_{L^2}^2+\frac{\delta^\frac12\eta}{2}\|\nabla H_t\|_{L^2}^2
+C\|\theta\|_{L^\infty}^2\|\nabla H\|_{L^2}^2+\delta^\frac32\|\nabla^2H\|_{L^2}^4+C\|\nabla H\|_{L^2}^4.
\end{align}
Combining $\eqref{f19}$, $\eqref{f20}$, and $\eqref{f28}$ multiplied by $\delta$ altogether, we thus deduce that, for all $0<\delta\le \min\Big\{1,
\frac{\eta^2}{4C_2}, \frac{c_v^2}{16}\Big\}$,
\begin{align}\label{f29}
&\frac{d}{dt}\int\Big(\frac{\mu}{2}|\nabla u|^2+\frac{\kappa\delta}{2}|\nabla\theta|^2+\eta|\nabla H|^2\Big)dx
+\int\Big(\rho|\dot{u}|^2+|H_t|^2+\frac{\eta^2}{2}|\nabla^2H|^2
+\frac{c_v\delta}{2}\rho|\dot{\theta}|^2
\Big)dx\nonumber\\
&\le \frac{d}{dt}\int\Big(R\rho\theta\divv u+\frac12|H|^2\divv u-H\cdot\nabla u\cdot H
+\lambda\theta\delta(\divv u)^2+\frac{\mu}{2}\delta\theta|\nabla u+(\nabla u)^{tr}|^2\Big)dx\nonumber\\
&\quad+\frac{d}{dt}\int\delta\theta|\curl H|^2dx
+\delta^{-\frac32}\big(\|\theta\nabla u\|_{L^2}^2+\|\nabla u\|_{L^2}^2+\|\nabla u\|_{L^2}^4\|\nabla\theta\|_{L^2}^2\big)
+2\delta^2\|\nabla\dot{u}\|_{L^2}^2\nonumber\\
&\quad+6\delta^\frac32\|\nabla u\|_{L^4}^4+C\|\nabla H\|_{L^2}^4
+C\|\nabla u\|_{L^2}^4\|\nabla H\|_{L^2}^2+C\|H\|_{L^2}\|\nabla H\|_{L^2}^3+\frac{\delta^\frac32\eta}{2}\|\nabla H_t\|_{L^2}^2\nonumber\\[3pt]
&\quad+\delta^\frac52\|\nabla^2H\|_{L^2}^4
+C\delta\|\theta\|_{L^\infty}^2\|\nabla H\|_{L^2}^2.
\end{align}

4. Operating $\partial_t+{\rm div}(u\cdot)$ to the $j$-th component of $\eqref{a1}_2$ and
multiplying the resulting equation by $\dot{u}^j$ , one gets by some calculations that
\begin{align}\label{f30}
\frac{1}{2}\frac{d}{dt}\int\rho|\dot{u}|^2dx
&=-\int\dot{u}^j(\partial_jp_t+{\rm div}(u\partial_jp))dx+\mu\int\dot{u}^j(\partial_t\Delta u^j+\divv(u\Delta u^j))dx
\nonumber\\
&\quad+(\mu+\lambda)\int\dot{u}^j(\partial_t\partial_j(\divv u)
+\divv(u\partial_j(\divv u)))dx\nonumber\\
&\quad-\frac12\int\dot{u}^j\big(\partial_t\partial_j|H|^2+{\rm div}(u\partial_j|H|^2)\big)dx\nonumber\\
&\quad+\int\dot{u}^j\big[\partial_t\partial_i(H^iH^j)
+\divv\big(u\partial_i(H^iH^j)\big)\big]dx\triangleq\sum_{i=1}^5K_i.
\end{align}
It follows from integration by parts, $\eqref{a1}_1$, and Young's inequality that
\begin{align}\label{z3.32}
K_1&=\int\big(\partial_j\dot{u}^jp_t+\partial_j pu\cdot\nabla\dot{u}^j\big)dx\nonumber\\
&=\int\partial_j\dot{u}^jp_tdx-\int p\partial_j\big(u\cdot\nabla\dot{u}^j\big)dx\nonumber\\
&=\int\partial_j\dot{u}^j\big[(\rho\theta)_t+\divv(\rho u\theta)\big]dx
-\int\rho\theta\partial_ju\cdot\nabla\dot{u}^jdx\nonumber\\
&=\int\rho\dot{\theta}\partial_j\dot{u}^jdx
-\int\rho\theta\partial_ju\cdot\nabla\dot{u}^jdx\nonumber\\
&\le \frac{\mu}{4}\|\nabla\dot{u}\|_{L^2}^2+C\|\theta\nabla u\|_{L^2}^2+C\|\sqrt{\rho}\dot{\theta}\|_{L^2}^2.
\end{align}
Integration by parts leads to
\begin{align}
K_2&=-\mu\int\big(\partial_i\dot{u}^j\partial_t\partial_iu^j+\Delta u^ju\cdot\nabla\dot{u}^j\big)dx\nonumber\\
&=-\mu\int\big(|\nabla\dot{u}|^2-\partial_i\dot{u}^ju^k\partial_k\partial_iu^j-\partial_i\dot{u}^j\partial_iu^k\partial_ku^j
+\Delta u^ju\cdot\nabla\dot{u}^j\big)dx\nonumber\\
&=-\mu\int\big(|\nabla\dot{u}|^2+\partial_i\dot{u}^j\partial_ku^k\partial_iu^j-\partial_i\dot{u}^j\partial_iu^k\partial_ku^j
-\partial_iu^j\partial_iu^k\partial_k\dot{u}^j\big)dx\nonumber\\
&\le -\frac{3\mu}{4}\|\nabla\dot{u}\|_{L^2}^2+C\|\nabla u\|_{L^4}^4.
\end{align}
Similarly, one has
\begin{align}
K_3&\le -\frac{\mu+\lambda}{2}\|\divv\dot{u}\|_{L^2}^2+C\|\nabla u\|_{L^4}^4.
\end{align}
Integrating by parts together with H\"older's
inequality and Sobolev's inequality leads to
\begin{align}
K_4&= -\frac12\int\dot{u}^j\big[\partial_t\partial_j|H|^2
+\divv\big(u\partial_j|H|^2\big)\big]dx\nonumber\\
&=\frac12\int\big(\partial_j\dot{u}^j\partial_t|H|^2
+\partial_i\dot{u}^ju^i\partial_j|H|^2\big)dx\nonumber\\
&\le C\|\nabla\dot{u}\|_{L^2}\big(\|H\|_{L^6}\|H_t\|_{L^3}+\|u\|_{L^6}\|H\|_{L^6}\|\nabla H\|_{L^6}\big)\nonumber\\
&\le \frac{\mu}{16}\|\nabla\dot{u}\|_{L^2}^2+C\|H\|_{L^6}^2\|H_t\|_{L^3}^2
+C\|u\|_{L^6}^2\|H\|_{L^6}^2\|\nabla H\|_{L^6}^2\nonumber\\
&\le \frac{\mu}{16}\|\nabla\dot{u}\|_{L^2}^2+C\|\nabla H\|_{L^2}^2\|H_t\|_{L^2}\|\nabla H_t\|_{L^2}+C\|\nabla u\|_{L^2}^2\|\nabla H\|_{L^2}^2\|\nabla^2H\|_{L^2}^2\nonumber\\
&\le \frac{\mu}{16}\|\nabla\dot{u}\|_{L^2}^2+\delta\|\nabla H_t\|_{L^2}^2
+C\|\nabla H\|_{L^2}^4\|H_t\|_{L^2}^2+C\|\nabla^2H\|_{L^2}^2.
\end{align}
Similarly, we obtain that
\begin{align}\label{z3.36}
K_5\le \frac{\mu}{16}\|\nabla\dot{u}\|_{L^2}^2+\delta\|\nabla H_t\|_{L^2}^2
+C\|\nabla H\|_{L^2}^4\|H_t\|_{L^2}^2+C\|\nabla^2H\|_{L^2}^2.
\end{align}
Substituting \eqref{z3.32}--\eqref{z3.36} into \eqref{f30}, we arrive at
\begin{align}\label{f35}
&\frac{1}{2}\frac{d}{dt}\int\rho|\dot{u}|^2dx+\frac{\mu}{2}\int|\nabla\dot{u}|^2dx\nonumber\\
&\le 2\delta\|\nabla H_t\|_{L^2}^2+C\|\theta\nabla u\|_{L^2}^2+C\|\sqrt{\rho}\dot{\theta}\|_{L^2}^2
+C\|\nabla H\|_{L^2}^4\|H_t\|_{L^2}^2+C\|\nabla^2H\|_{L^2}^2.
\end{align}

5. Differentiating \eqref{a1}$_4$ with respect to $t$ yields that
\begin{align}\label{w36}
H_{tt}-\eta\Delta H_t=\big(H\cdot\nabla u-u\cdot\nabla H-H{\rm div}u\big)_t.
\end{align}
Multiplying \eqref{w36} by $H_t$ and noting that $u_t=\dot{u}-u\cdot\nabla u$, we obtain from integration by parts, H\"older's inequality, Sobolev's inequality, Gagliardo-Nirenberg inequality, Young's inequality, \eqref{lz3.3}, and \eqref{ee1} that
\begin{align}\label{rrf}
&\frac12\frac{d}{dt}\int|H_t|^2dx+\eta\int|\nabla H_t|^2dx\nonumber\\
&=\int\big(H_t\cdot\nabla u\cdot H_t+H\cdot\nabla u_t\cdot H_t-u_t\cdot\nabla H\cdot H_t\big)dx\nonumber\\
&\quad-\int\big(u\cdot\nabla H_t
\cdot H_t-|H_t|^2\divv u-H\cdot H_t\divv u_t\big)dx\nonumber\\
&=\int\big(H_t\cdot\nabla u\cdot H_t-H\cdot\nabla H_t\cdot u_t-u\cdot\nabla H_t\cdot H_t
+u_t\cdot\nabla H_t\cdot H-|H_t|^2\divv u\big)dx\nonumber\\
&=\int\big(H_t\cdot\nabla u\cdot H_t-H\cdot\nabla H_t\cdot\dot{u}
+H\cdot\nabla H_t\cdot(u\cdot\nabla u)\big)dx\nonumber\\
&\quad+\int\big(-u\cdot\nabla H_t\cdot H_t
+\dot{u}\cdot\nabla H_t\cdot H-(u\cdot\nabla u)\cdot\nabla H_t\cdot H-|H_t|^2
\divv u\big)dx\nonumber\\
&\le C\int|H_t|^2|\nabla u|dx+C\int|u||H||\nabla u||\nabla H_t|dx\nonumber\\
&\quad+C\int|H||\nabla H_t||\dot{u}|dx+C\int|u||\nabla H_t||H_t|dx\nonumber\\
&\le C\|\nabla u\|_{L^2}\|H_t\|_{L^4}^2+C\|\nabla H_t\|_{L^2}
\|H\|_{L^3}\|\dot{u}\|_{L^6}+C\|\nabla H_t\|_{L^2}\|u\|_{L^6}\|H_t\|_{L^3}
\nonumber\\
&\quad
+C\|H\|_{L^{12}}\|\nabla H_t\|_{L^2}\|u\|_{L^6}\|\nabla u\|_{L^4}\nonumber\\
&\le C\|\nabla u\|_{L^2}\|H_t\|_{L^2}^\frac12\|\nabla H_t\|_{L^2}^\frac32
+C\|\nabla H_t\|_{L^2}\|H\|_{L^2}^\frac12\|\nabla H\|_{L^2}^\frac12\|\nabla\dot{u}\|_{L^2}\nonumber\\
&\quad+C\|H\|_{L^2}^\frac14\|\nabla H\|_{L^2}^\frac14\|\nabla^2H\|_{L^2}^\frac12\|\nabla H_t\|_{L^2}\|\nabla u\|_{L^2}\|\nabla u\|_{L^4}\nonumber\\
&\le C\|\nabla u\|_{L^2}\|H_t\|_{L^2}^\frac12\|\nabla H_t\|_{L^2}^\frac32
+C\|\nabla H_t\|_{L^2}\|\nabla\dot{u}\|_{L^2}
+C\|\nabla^2H\|_{L^2}^\frac12\|\nabla H_t\|_{L^2}\|\nabla u\|_{L^4}
\nonumber\\
&\le \frac{\eta}{2}\|\nabla H_t\|_{L^2}^2+c_2\|\nabla\dot{u}\|_{L^2}^2
+C\|\nabla u\|_{L^4}^4+C\|\nabla^2H\|_{L^2}^2+C\|\nabla u\|_{L^2}^4\|H_t\|_{L^2}^2,
\end{align}
which implies that
\begin{align}\label{f37}
\frac{d}{dt}\int|H_t|^2dx+\eta\int|\nabla H_t|^2dx
\le c_2\|\nabla\dot{u}\|_{L^2}^2
+C\|\nabla u\|_{L^4}^4+C\|\nabla^2H\|_{L^2}^2+C\|\nabla u\|_{L^2}^4\|H_t\|_{L^2}^2.
\end{align}
Adding \eqref{f35} multiplied by $\big(\frac{2c_2}{\mu}+1\big)$ to \eqref{f37}, we get that
\begin{align}\label{f38}
&\frac{d}{dt}\int\Big[\Big(\frac{c_2}{\mu}+\frac12\Big)\rho|\dot{u}|^2
+|H_t|^2\Big]dx
+\int\Big(\frac{\mu}{2}|\nabla \dot{u}|^2+\frac{\eta}{2}|\nabla H_t|^2\Big)dx\nonumber\\
&\le C\|\nabla u\|_{L^4}^4+C\|\nabla^2H\|_{L^2}^2
+C\|\nabla u\|_{L^2}^4\|H_t\|_{L^2}^2
+C\|\theta\nabla u\|_{L^2}^2+C\|\sqrt{\rho}\dot{\theta}\|_{L^2}^2\nonumber\\
&\le c_3\big(\|\nabla u\|_{L^4}^4+\|\nabla^2H\|_{L^2}^2+\|H_t\|_{L^2}^2+\|\sqrt{\rho}\dot{\theta}\|_{L^2}^2
+\|\theta\nabla u\|_{L^2}^2\big).
\end{align}
Thus, combining \eqref{f29} and \eqref{f38} multiplied by $2\delta^\frac{3}{2}$ together, we conclude \eqref{f12} as long as $\delta$ is chosen suitably small such that
\begin{align*}
0<\delta\le\min\left\{\frac{\eta^2}{4C_2}, \frac{c_v^2}{16}, \frac{1}{\sqrt[3]{4}}, \frac{1}{2C_2},
\frac{c_v^2}{(8c_3)^2}, \frac{1}{(4c_3)^\frac32}, \frac{\mu\eta}{4(2c_2+\mu)}, \frac{\eta^\frac43}{(8c_3)^\frac23}\right\}.
\end{align*}
The proof of Lemma \ref{l32} is complete.  \hfill $\Box$

Now we define
\begin{align}\label{f39}
\sigma(t)\triangleq\min\{t_0, t\},
\end{align}
where $t_0\in (0, T_*)$ is the same as in Lemma \ref{l31}.
\begin{lemma}\label{l33}
Let $\sigma(t)$ be as in \eqref{f39} and let $(\rho, u, \theta, H)$ be the smooth solution of \eqref{a1}--\eqref{a3} satisfying \eqref{ee1}, then it holds that
\begin{align}\label{f42}
&\sup_{0\le t\le T}\big[\sigma(t)\big(\|\nabla u\|_{L^2}^2
+\|\nabla\theta\|_{L^2}^2+\|\nabla H\|_{L^2}^2+\|\sqrt{\rho}\dot{u}\|_{L^2}^2
+\|H_t\|_{L^2}^2\big)\big]\nonumber\\
&\quad+\int_0^T\sigma(t)\big(\|\sqrt{\rho}\dot{u}\|_{L^2}^2
+\|\sqrt{\rho}\dot{\theta}\|_{L^2}^2
+\|\nabla\dot{u}\|_{L^2}^2+\|H_{tt}\|_{L^2}^2+\|H_t\|_{L^2}^2+\|\nabla^2H\|_{L^2}^2+\|\nabla H_t\|_{L^2}^2\big)dt\nonumber\\
&\le C\Big(t_0+\mathbb{E}_0^\frac{4}{5}\Big).
\end{align}
\end{lemma}
{\it Proof.}
1. Multiplying \eqref{f12} by $\sigma(t)$, and integrating the resultant with respect to $t$, we have
\begin{align}\label{f43}
&\sup_{0\le t\le T}\big[\sigma\big(\|\nabla u\|_{L^2}^2
+\delta\|\nabla\theta\|_{L^2}^2+\|\nabla H\|_{L^2}^2+\delta^\frac32\|\sqrt{\rho}\dot{u}\|_{L^2}^2
+\delta^\frac32\|H_t\|_{L^2}^2\big)\big]\nonumber\\
&\quad+\int_0^T\sigma\big(\|\sqrt{\rho}\dot{u}\|_{L^2}^2+\|H_t\|_{L^2}^2
+\|\nabla^2H\|_{L^2}^2+\delta\|\sqrt{\rho}\dot{\theta}\|_{L^2}^2
+\delta^\frac32\|\nabla\dot{u}\|_{L^2}^2+\delta^\frac32\|\nabla H_t\|_{L^2}^2\big)dt\nonumber\\
&\le C\int_0^{t_0}\big(\|\nabla u\|_{L^2}^2
+\delta\|\nabla\theta\|_{L^2}^2+\|\nabla H\|_{L^2}^2+\delta^\frac32\|\sqrt{\rho}\dot{u}\|_{L^2}^2
+\delta^\frac{3}{2}\|H_t\|_{L^2}^2\big)dt\nonumber\\
&\quad+C\int_0^{t_0}\int\big(\rho\theta|\divv u|+\theta|\nabla u|^2+|\nabla u||H|^2+\theta|\nabla H|^2\big)dxdt\nonumber\\
&\quad+C\sigma\int\big(\rho\theta|\nabla u|+|\nabla u||H|^2+\delta\theta|\nabla u|^2+\delta\theta|\nabla H|^2\big)dx\nonumber\\
&\quad+C\int_0^T\sigma\big(\|\theta\nabla u\|_{L^2}^2+\|\nabla u\|_{L^2}^2
+\|\nabla u\|_{L^2}^4\|\nabla\theta\|_{L^2}^2\big)dt+C\delta^\frac52\int_0^T\sigma\|\nabla^2H\|_{L^2}^4dt\nonumber\\
&\quad+C\delta\int_0^T\sigma\|\theta\|_{L^\infty}^2\|\nabla H\|_{L^2}^2dt
+C\delta^\frac32\int_0^T\sigma\|\nabla u\|_{L^4}^4dt\nonumber\\
&\quad+C\int_0^T\sigma\big(\|\nabla H\|_{L^2}^2+\|\nabla u\|_{L^2}^4
+\|H\|_{L^2}^2\big)\|\nabla H\|_{L^2}^2dt.
\end{align}
Let us deal with the right-hand side terms in \eqref{f43} as follows. By virtue of \eqref{w1}, one has
\begin{align}\label{f44}
&C\int_0^{t_0}\big(\|\nabla u\|_{L^2}^2
+\delta\|\nabla\theta\|_{L^2}^2+\|\nabla H\|_{L^2}^2+\delta^\frac{3}{2}\|\sqrt{\rho}\dot{u}\|_{L^2}^2
+\delta^\frac{3}{2}\|H_t\|_{L^2}^2\big)dt\nonumber\\
&\le Ct_0\sup_{0\le t\le t_0}\big(\|\nabla u\|_{L^2}^2
+\delta\|\nabla\theta\|_{L^2}^2+\|\nabla H\|_{L^2}^2+\delta^\frac{3}{2}\|\sqrt{\rho}\dot{u}\|_{L^2}^2
+\delta^\frac{3}{2}\|H_t\|_{L^2}^2\big)\nonumber\\
&\le C(M_*)t_0.
\end{align}
It follows from \eqref{w1}, \eqref{ee1}, H\"older's, Young's, and Gagliardo-Nirenberg inequalities that
\begin{align}\label{f45}
&C\int_0^{t_0}\int\big(\rho\theta|{\rm div}u|+\theta|\nabla u|^2+|\nabla u||H|^2+\theta|\nabla H|^2\big)dxdt\nonumber\\
&\le C(\bar{\rho})\int_0^{t_0}\|\nabla u\|_{L^2}\big(\|\sqrt{\rho}\theta\|_{L^2}
+\|H\|_{L^4}^2\big)dt+C\int_0^{t_0}\|\theta\|_{L^\infty}\big(\|\nabla u\|_{L^2}^2+\|\nabla H\|_{L^2}^2\big)dt\nonumber\\
&\le C\int_0^{t_0}\|\nabla u\|_{L^2}\Big(\|\sqrt{\rho}\theta\|_{L^2}
+\|H\|_{L^2}^\frac{1}{2}\|\nabla H\|_{L^2}^\frac{3}{2}\Big)dt\nonumber\\
&\quad+C\int_0^{t_0}\|\nabla \theta\|_{L^2}^\frac{1}{2}\|\nabla^2\theta\|_{L^2}^\frac{1}{2}\big(
\|\nabla u\|_{L^2}^2+\|\nabla H\|_{L^2}^2\big)dt\nonumber\\
&\le C\int_0^{t_0}\big(\|\nabla u\|_{L^2}^2+\|\sqrt{\rho}\theta\|_{L^2}^2+\|\nabla H\|_{L^2}^2
+\|H\|_{L^2}^2\|\nabla u\|_{L^2}^4\big)dt\nonumber\\
&\quad+C\int_0^{t_0}\big(\|\nabla\theta\|_{L^2}^2+\|\nabla^2\theta\|_{L^2}^2
+\|\nabla u\|_{L^2}^4+\|\nabla H\|_{L^2}^4\big)dt\nonumber\\
&\le C(\bar{\rho}, M_*)t_0.
\end{align}
For the third term in the right hand side of \eqref{f43}, we need the following fact
\begin{align}\label{f46}
C\sigma\int|\nabla u||H|^2dx&\le C\sigma\|\nabla u\|_{L^2}\|H\|_{L^4}^2\nonumber\\
&\le \frac{1}{4}\sigma\|\nabla u\|_{L^2}^2+C\sigma\|H\|_{L^2}\|\nabla H\|_{L^2}^3\nonumber\\
&\le \frac{1}{4}\sigma\|\nabla u\|_{L^2}^2+C\delta^{-1}\sigma\|H\|_{L^2}^2\|\nabla H\|_{L^2}^2+\delta\sigma\|\nabla H\|_{L^2}^4\nonumber\\
&\le \frac{1}{4}\sigma\|\nabla u\|_{L^2}^2+C(M_*)\Big(\mathbb{E}_0^\frac23
+\delta\Big)\sup_{0\le t\le T}\big(\sigma\|\nabla H\|_{L^2}^2\big).
\end{align}
By \eqref{w7}, \eqref{ee1}, and Gagliardo-Nirenberg inequality, we have
\begin{align}\label{f47}
\|\nabla u\|_{L^6}
&\le C\big(\|\sqrt{\rho}\dot{u}\|_{L^2}+\|\nabla\theta\|_{L^2}+\|H\|_{L^\infty}\|\nabla H\|_{L^2}\big)\nonumber\\
&\le C\big(\|\sqrt{\rho}\dot{u}\|_{L^2}+\|\nabla\theta\|_{L^2}+\|\nabla H\|_{L^2}^\frac{3}{2}\|\nabla^2H\|_{L^2}^\frac{1}{2}\big)\nonumber\\
&\le C\big(\|\sqrt{\rho}\dot{u}\|_{L^2}+\|\nabla\theta\|_{L^2}
+\|\nabla^2H\|_{L^2}+\|\nabla H\|_{L^2}^3\big).
\end{align}
Moreover, it follows from $\eqref{a1}_4$, the standard $L^2$-estimate of elliptic equations, and Lemma \ref{l22} that
\begin{align*}
\|\nabla^2H\|_{L^2}&\le C\|H_t-H\cdot\nabla u+u\cdot\nabla H+H\divv u\|_{L^2}\nonumber\\
&\le C\|H_t\|_{L^2}+C\|u\|_{L^6}\|\nabla H\|_{L^3}+C\|H\|_{L^\infty}\|\nabla u\|_{L^2}\nonumber\\
&\le C\|H_t\|_{L^2}+C\|\nabla u\|_{L^2}\|\nabla H\|_{L^2}^\frac12\|\nabla^2H\|_{L^2}^\frac12\nonumber\\
&\le \frac12\|\nabla^2H\|_{L^2}+C\|H_t\|_{L^2}+C\|\nabla u\|_{L^2}^2\|\nabla H\|_{L^2},
\end{align*}
which leads to
\begin{align}\label{f49}
\|\nabla^2H\|_{L^2}\le C\|H_t\|_{L^2}+C\|\nabla u\|_{L^2}^2\|\nabla H\|_{L^2}.
\end{align}
Inserting \eqref{f49} into \eqref{f47} gives that
\begin{align}\label{f50}
\|\nabla u\|_{L^6}\le C\|\sqrt{\rho}\dot{u}\|_{L^2}+C\|H_t\|_{L^2}+C\|\nabla\theta\|_{L^2}
+C\|\nabla H\|_{L^2}^3+C\|\nabla u\|_{L^2}^2\|\nabla H\|_{L^2}.
\end{align}
Thus, in terms of \eqref{f50}, we estimate
\begin{align*}
C\int\theta|\nabla u|^2dx
&\le C\delta\|\theta\|_{L^6}\|\nabla u\|_{L^2}^\frac{3}{2}\|\nabla u\|_{L^6}^\frac12\nonumber\\
&\le C\big(\|\sqrt{\rho}\dot{u}\|_{L^2}+\|H_t\|_{L^2}+\|\nabla\theta\|_{L^2}
+\|\nabla H\|_{L^2}^3+\|\nabla u\|_{L^2}^2\|\nabla H\|_{L^2}\big)^\frac12\nonumber\\
&\quad\times\|\nabla u\|_{L^2}^\frac{1}{2}\|\nabla\theta\|_{L^2}\|\nabla u\|_{L^2}\nonumber\\
&\le C\delta^{-\frac78}\|\nabla u\|_{L^2}^2\|\nabla\theta\|_{L^2}^2
+C\delta^{-\frac18}\|\nabla u\|_{L^2}^2+\delta^{\frac{15}{8}}\|\sqrt{\rho}\dot{u}\|_{L^2}^2
+\delta^{\frac{15}{8}}\|H_t\|_{L^2}^2\nonumber\\
&\quad+\delta^{\frac{15}{8}}\|\nabla\theta\|_{L^2}^2
+\delta^{\frac{15}{8}}\big(\|\nabla u\|_{L^2}^4+\|\nabla H\|_{L^2}^4\big)\|\nabla H\|_{L^2}^2,
\end{align*}
that is
\begin{align}\label{f51}
C\delta\sigma\int\theta|\nabla u|^2dx&\le C\delta^\frac{1}{8}\sup_{0\le t\le T}\big(\sigma\|\nabla u\|_{L^2}^2\big)
+\delta^\frac{15}{8}\sup_{0\le t\le T}\big(\delta\sigma\|\nabla\theta\|_{L^2}^2\big)\nonumber\\
&\quad+\delta^\frac{11}{8}\sup_{0\le t\le T}\big(\delta^\frac{3}{2}\sigma\|\sqrt{\rho}\dot{u}\|_{L^2}^2\big)
+\delta^\frac{11}{8}\sup_{0\le t\le T}\big(\delta^\frac{3}{2}\sigma\|H_t\|_{L^2}^2\big)\nonumber\\
&\quad+C(M_*)\delta^\frac{15}{8}\sup_{0\le t\le T}\big(\sigma\|\nabla H\|_{L^2}^2\big).
\end{align}
By Sobolev's inequality and H\"older's inequality, we have
\begin{align*}
C\int\theta|\nabla H|^2dx
&\le C\|\theta\|_{L^6}\|\nabla H\|_{L^2}\|\nabla H\|_{L^3}\nonumber\\
&\le C\|\nabla\theta\|_{L^2}\|\nabla H\|_{L^2}\|\nabla H\|_{L^2}^\frac{1}{2}\|\nabla^2H\|_{L^2}^\frac{1}{2}
\nonumber\\
&\le C\delta^{-\frac{7}{8}}\|\nabla\theta\|_{L^2}^2\|\nabla H\|_{L^2}^2
+\delta^\frac{7}{8}\|\nabla^2H\|_{L^2}^2+C\delta^\frac{7}{8}\|\nabla H\|_{L^2}^2\nonumber\\
&\le C\delta^{-\frac{7}{8}}\|\nabla H\|_{L^2}^2
+C\delta^\frac{7}{8}\|H_t\|_{L^2}^2+C(M_*)\delta^\frac{7}{8}\|\nabla H\|_{L^2}^2,
\end{align*}
which yields that
\begin{align}\label{f52}
C\delta\sigma\int\theta|\nabla H|^2dx&\le C\delta^\frac{1}{8}\sup_{0\le t\le T}\big(\sigma\|\nabla H\|_{L^2}^2\big)
+\delta^\frac{3}{8}\sup_{0\le t\le T}\big(\sigma\delta^\frac{3}{2}\|H_t\|_{L^2}^2\big).
\end{align}
Thus, combining \eqref{f46}, \eqref{f51}, and \eqref{f52} altogether, we derive that
\begin{align}\label{f53}
&C\sigma\int\big(\rho\theta|\nabla u|+|\nabla u||H|^2+\delta\theta|\nabla u|^2+\delta\theta|\nabla H|^2\big)dx\nonumber\\
&\le C\sup_{0\le t\le T}\|\sqrt{\rho}\theta\|_{L^2}^2
+\Big(\frac{1}{4}+C\delta^\frac{1}{8}\Big)\sup_{0\le t\le T}\big(\sigma\|\nabla u\|_{L^2}^2\big)\nonumber\\
&\quad+C\Big(\mathbb{E}_0^\frac23+\delta^\frac{1}{8}\Big)\sup_{0\le t\le T}\big(\sigma\|\nabla H\|_{L^2}^2\big)
+\delta^\frac{3}{8}\sup_{0\le t\le T}\big(\sigma\delta^\frac{3}{2}\|H_t\|_{L^2}^2\big)
+\delta^\frac{15}{8}\sup_{0\le t\le T}\big(\sigma\delta\|\nabla\theta\|_{L^2}^2\big)\nonumber\\
&\quad+\delta^\frac{11}{8}\sup_{0\le t\le T}\big(\sigma\delta^\frac{3}{2}\|\sqrt{\rho}\dot{u}\|_{L^2}^2\big)\nonumber\\
&\le C\Big(t_0+\mathbb{E}_0^\frac{4}{5}\Big)+\Big[\frac{1}{4}
+C_3\Big(\mathbb{E}_0^\frac{2}{3}+\delta^\frac{1}{8}\Big)\Big]
\sup_{0\le t\le T}\big(\sigma\|\nabla u\|_{L^2}^2\big)
+\delta^\frac{11}{8}\sup_{0\le t\le T}\big(\sigma\delta^\frac{3}{2}\|\sqrt{\rho}\dot{u}\|_{L^2}^2\big)\nonumber\\
&\quad+C_3\Big(\mathbb{E}_0^\frac{2}{3}+\delta^\frac{1}{8}\Big)\sup_{0\le t\le T}\big(\sigma\|\nabla H\|_{L^2}^2\big)
+C_3\Big(\mathbb{E}_0^\frac{2}{3}\delta^{-1}+\delta^\frac{15}{8}\Big)\sup_{0\le t\le T}\big(\sigma\delta\|\nabla\theta\|_{L^2}^2\big)\nonumber\\
&\quad+\delta^\frac{3}{8}\sup_{0\le t\le T}\big(\sigma\delta^\frac{3}{2}\|H_t\|_{L^2}^2\big).
\end{align}
Let
\begin{align*}
\delta=\delta_1\triangleq\min\left\{\frac{\eta^2}{4C_2}, \frac{c_v^2}{16}, \frac{1}{\sqrt[3]{4}}, \frac{1}{2C_2},
\frac{c_v^2}{(8c_3)^2}, \frac{1}{(4c_3)^\frac32}, \frac{\mu\eta}{4(2c_2+\mu)}, \frac{\eta^\frac43}{(8c_3)^\frac23}, \frac{1}{(8C_3)^8}\right\},
\end{align*}
substituting \eqref{f44}, \eqref{f45}, and \eqref{f53} into \eqref{f43}, we conclude that, if $\mathbb{E}_0\le \varepsilon_2\triangleq\min\left\{\varepsilon_1,
\Big(\frac{\delta_1}{8C_3}\Big)^\frac32\right\}$,
\begin{align}\label{f54}
&\sup_{0\le t\le T}\big[\sigma\big(\|\nabla u\|_{L^2}^2
+\|\nabla\theta\|_{L^2}^2+\|\nabla H\|_{L^2}^2+\|\sqrt{\rho}\dot{u}\|_{L^2}^2
+\|H_t\|_{L^2}^2\big)\big]\nonumber\\
&\quad+\int_0^T\sigma\big(\|\sqrt{\rho}\dot{u}\|_{L^2}^2+\|H_t\|_{L^2}^2
+\|\nabla^2H\|_{L^2}^2+\|\sqrt{\rho}\dot{\theta}\|_{L^2}^2
+\|\nabla\dot{u}\|_{L^2}^2+\|\nabla H_t\|_{L^2}^2\big)dt\nonumber\\
&\le C\Big(t_0+\mathbb{E}_0^\frac{4}{5}\Big)+C\int_0^T\sigma\big(\|\theta\nabla u\|_{L^2}^2+\|\nabla u\|_{L^2}^2
+\|\nabla u\|_{L^2}^4\|\nabla\theta\|_{L^2}^2\big)dt\nonumber\\
&\quad+C\delta\int_0^T\sigma\|\theta\|_{L^\infty}^2\|\nabla H\|_{L^2}^2dt
+C\delta^\frac52\int_0^T\sigma\|\nabla^2H\|_{L^2}^4dt+C\delta^\frac32\int_0^T\sigma\|\nabla u\|_{L^4}^4dt\nonumber\\
&\quad+C\int_0^T\sigma\big(\|\nabla H\|_{L^2}^2+\|\nabla u\|_{L^2}^4
+\|H\|_{L^2}^2\big)\|\nabla H\|_{L^2}^2dt\nonumber\\
&\triangleq C\Big(t_0+\mathbb{E}_0^\frac{4}{5}\Big)+\sum_{i=1}^5\bar{R}_i.
\end{align}
It follows from H\"older's, Young's, Gagliardo-Nirenberg inequalities, \eqref{ee1}, \eqref{g12}, and \eqref{f5} that
\begin{align*}
\bar{R}_1&\le \sup_{0\le t\le T}\big(\sigma\|\nabla u\|_{L^2}^2\big)\int_0^T\big(\|\nabla \theta\|_{L^2}^2
+\|\nabla^2\theta\|_{L^2}^2\big)dt+C(M_*)\int_0^T\big(\|\nabla u\|_{L^2}^2+\|\nabla\theta\|_{L^2}^2\big)dt\nonumber\\
&\le C\mathbb{E}_0^\frac{2}{3}\sup_{t_0\le t\le T}\big(\|\nabla u\|_{L^2}^2+\|\nabla\theta\|_{L^2}^2
+\|\nabla H\|_{L^2}^2\big)+C\Big(t_0+\mathbb{E}_0^\frac{4}{5}\Big)
+C\mathbb{E}_0^\frac{2}{3}\sup_{0\le t\le T}\big(\sigma\|\nabla u\|_{L^2}^2\big)\nonumber\\
&\le  C\Big(t_0+\mathbb{E}_0^\frac{4}{5}\Big)+C\mathbb{E}_0^\frac{2}{3}\sup_{0\le t\le T}\sigma\big(\|\nabla u\|_{L^2}^2
+\|\nabla\theta\|_{L^2}^2+\|\nabla H\|_{L^2}^2\big),\\
\bar{R}_2&\le C(M_*)\sup_{0\le t\le T}\big(\sigma\|\nabla H\|_{L^2}^2\big)\int_0^T\big(\|\nabla\theta\|_{L^2}^2
+\|\nabla^2\theta\|_{L^2}^2\big)dt\nonumber\\
\bar{R}_3&\le C\sup_{0\le t\le T}\big(\sigma\|\nabla^2H\|_{L^2}^2\big)\int_0^T\|\nabla^2H\|_{L^2}^2dt\nonumber\\
&\le C\mathbb{E}_0^\frac32\sup_{0\le t\le T}\sigma\big(\|H_t\|_{L^2}^2+\|\nabla u\|_{L^2}^2+\|\nabla H\|_{L^2}^2\big),\\
\bar{R}_4&\le C(M_*)\int_0^T\sigma\|\nabla u\|_{L^2}\big(\|\sqrt{\rho}\dot{u}\|_{L^2}+\|H_t\|_{L^2}+\|\nabla\theta\|_{L^2}
+\|\nabla H\|_{L^2}\big)^3dt\nonumber\\
&\le \sup_{0\le t\le T}\big(\sigma\|H_t\|_{L^2}^2\big)\int_0^T\|H_t\|_{L^2}^2dt
+\sup_{0\le t\le T}\big(\sigma\|\sqrt{\rho}\dot{u}\|_{L^2}^2\big)
\int_0^T\|\sqrt{\rho}\dot{u}\|_{L^2}^2dt\nonumber\\
&\quad+C\int_0^T\big(\|\nabla u\|_{L^2}^4
+\|\nabla\theta\|_{L^2}^4+\|\nabla H\|_{L^2}^4\big)dt
\nonumber\\
&\le C\Big(t_0+\mathbb{E}_0^\frac{4}{5}\Big)+C\mathbb{E}_0^\frac{2}{3}\sup_{t_0\le t\le T}\big(\|\nabla u\|_{L^2}^2
+\|\nabla H\|_{L^2}^2+\|\nabla\theta\|_{L^2}^2\big)\nonumber\\
&\quad+2\mathbb{E}_0^\frac{2}{3}\sup_{0\le t\le T}\big[\sigma\big(\|\sqrt{\rho}\dot{u}\|_{L^2}^2
+\|H_t\|_{L^2}^2\big)\big]\nonumber\\
&\le C\mathbb{E}_0^\frac{2}{3}\sup_{0\le t\le T}\big[\sigma\big(\|\nabla u\|_{L^2}^2
+\|\nabla H\|_{L^2}^2+\|\nabla\theta\|_{L^2}^2+\|\sqrt{\rho}\dot{u}\|_{L^2}^2
+\|H_t\|_{L^2}^2\big)\big]+C\Big(t_0+\mathbb{E}_0^\frac{4}{5}\Big),\\
\bar{R}_5&\le C\int_0^T\sigma\big(\|\nabla u\|_{L^2}^4+\|\nabla H\|_{L^2}^4+\|H\|_{L^2}^2\|\nabla H\|_{L^2}^2\big)dt\nonumber\\
&\le C(M_*)\int_0^T\big(\|\nabla u\|_{L^2}^2+\|\nabla H\|_{L^2}^2\big)dt\nonumber\\
&\le C\Big(t_0+\mathbb{E}_0^\frac{4}{5}\Big)+C\mathbb{E}_0^\frac{2}{3}\sup_{t_0\le t\le T}\big(\|\nabla u\|_{L^2}^2
+\|\nabla\theta\|_{L^2}^2+\|\nabla H\|_{L^2}^2\big)\nonumber\\
&\le C\Big(t_0+\mathbb{E}_0^\frac{4}{5}\Big)+C\mathbb{E}_0^\frac{2}{3}
\sup_{0\le t\le T}\big[\sigma\big(\|\nabla u\|_{L^2}^2
+\|\nabla\theta\|_{L^2}^2+\|\nabla H\|_{L^2}^2\big)\big].
\end{align*}
Substituting $\bar{R}_1$--$\bar{R}_5$ into \eqref{f54}, one obtains that
\begin{align}
&\sup_{0\le t\le T}\big[\sigma\big(\|\nabla u\|_{L^2}^2
+\|\nabla\theta\|_{L^2}^2+\|\nabla H\|_{L^2}^2+\|\sqrt{\rho}\dot{u}\|_{L^2}^2
+\|H_t\|_{L^2}^2\big)\big]\nonumber\\
&\quad+\int_0^T\sigma\big(\|\sqrt{\rho}\dot{u}\|_{L^2}^2+\|H_t\|_{L^2}^2
+\|\nabla^2H\|_{L^2}^2+\|\sqrt{\rho}\dot{\theta}\|_{L^2}^2
+\|\nabla\dot{u}\|_{L^2}^2\big)dt+\int_0^T\sigma\|\nabla H_t\|_{L^2}^2dt\nonumber\\
&\le C_4\mathbb{E}_0^\frac{2}{3}\sup_{0\le t\le T}\big[\sigma\big(\|\nabla u\|_{L^2}^2
+\|\nabla H\|_{L^2}^2+\|\nabla\theta\|_{L^2}^2+\|\sqrt{\rho}\dot{u}\|_{L^2}^2
+\|H_t\|_{L^2}^2\big)\big]+C\Big(t_0+\mathbb{E}_0^\frac{4}{5}\Big).
\end{align}
If
\begin{align}
\mathbb{E}_0\le \varepsilon_3\triangleq\min
\left\{\varepsilon_2,\Big(\frac{1}{2C_4}\Big)^\frac32\right\},
\end{align}
 we then derive \eqref{f42}.
The proof of Lemma \ref{l33} is complete.  \hfill $\Box$

\begin{lemma}\label{l34}
Let $(\rho, u, \theta, H)$ be the smooth solution of \eqref{a1}--\eqref{a3} satisfying \eqref{ee1}, then it holds that
\begin{align}\label{3.49}
0\le \rho\le \frac{3\bar{\rho}}{2}.
\end{align}
\end{lemma}
{\it Proof.}
1. The first inequality of \eqref{3.49} is from \eqref{lz3.2}. Motivated by \cite{D1997} (see also \cite{WZ17,LZ20}), we next prove the second inequality of \eqref{3.49}. For any given $(x, t)\in \mathbb{R}^3\times[0, T]$, denote
\begin{align}
\rho^\delta(y,s)=\rho(y, s)+\delta\exp\left\{-\int_0^s\divv u(X(\tau; x, t), \tau)d\tau\right\}>0,
\end{align}
where
$X(s; x, t)$ is given by
\begin{align}
\left\{
\begin{array}{ll}
\displaystyle
\frac{d}{ds}X(s; x, t)=u(X(s; x, t), s), \quad 0\le s<t,\\
X(t; x, t)=x.
\end{array}
\right.
\end{align}
Using the fact that $\frac{d}{ds}(f(X(s; x, t), s)=(f_s+u\cdot\nabla f)(X(s; x, t), s)$,
 it follows from $\eqref{a1}_1$ that
 \begin{align}\label{tb1}
 \frac{d}{ds}\big(\ln(\rho^\delta(X(s; x, t), s)\big)=-\divv u(X(s; x, t), s),
 \end{align}
which leads to
\begin{align}
Y'(s)=g(s)+b'(s).
\end{align}
Here
\begin{align*}
&Y(s)=\ln\rho^\delta(X(s; x, t), s), \quad g(s)=-\frac{p(X(s; x, t), s)}{2\mu+\lambda},\\
&b(s)=-\frac{1}{2\mu+\lambda}\int_0^sF(X(\tau; x, t), \tau)d\tau,
\end{align*}
with $F$ as that of in \eqref{f3}.

2. Rewrite $\eqref{a1}_2$ as
\begin{align*}
\partial_t\big[\Delta^{-1}\divv(\rho u)\big]-(2\mu+\lambda)\divv u+p+\frac12|H|^2=-\Delta^{-1}\divv\divv(\rho u\otimes u)
-\Delta^{-1}\divv\divv(H\otimes H),
\end{align*}
which implies that
\begin{align}
F(X(\tau; x, t), \tau)
&=-\big[(-\Delta)^{-1}\divv(\rho u)\big]_{\tau}-(-\Delta)^{-1}\divv\divv(\rho u\otimes u)
-(-\Delta)^{-1}\divv\divv(H\otimes H)\nonumber\\
&=-\big[(-\Delta)^{-1}\divv(\rho u)\big]_{\tau}-u\cdot\nabla(-\Delta)^{-1}\divv(\rho u)
+u\cdot\nabla(-\Delta)^{-1}\divv(\rho u)\nonumber\\
&\quad-(-\Delta)^{-1}\divv\divv(\rho u\otimes u)
-(-\Delta)^{-1}\divv\divv(H\otimes H)\nonumber\\
&=-\frac{d}{d\tau}\big[(-\Delta)^{-1}\divv(\rho u)\big]+[u^i, R_{ij}](\rho u_j)-(-\Delta)^{-1}\divv\divv(H\otimes H),
\end{align}
where $[u^i, R_{ij}]=u^iR_{ij}-R_{ij}u^i$ and $R_{ij}=\partial_i(-\Delta)^{-1}\partial_j$ is the Riesz transform (we refer the reader to \cite[Chapter 3]{S70} for its definition and properties).
This yields that
\begin{align}\label{xyx}
b(t)-b(0)
&=\frac{1}{2\mu+\lambda}\int_0^t\Big[\frac{d}{d\tau}\big[(-\Delta)^{-1}\divv(\rho u)\big]-[u^i, R_{ij}](\rho u^j)
+(-\Delta)^{-1}\divv\divv(H\otimes H)\Big]d\tau\nonumber\\
&\quad+\frac{1}{2(2\mu+\lambda)}\int_0^t\|H\|_{L^\infty}^2d\tau\nonumber\\
&=\frac{1}{2\mu+\lambda}(-\Delta)^{-1}\divv(\rho u)-\frac{1}{2\mu+\lambda}(-\Delta)^{-1}\divv(\rho_0u_0)
-\frac{1}{2\mu+\lambda}\int_0^t[u^i, R_{ij}](\rho u^j)d\tau\nonumber\\
&\quad-\frac{1}{2\mu+\lambda}\int_0^t(-\Delta)^{-1}\divv\divv(H\otimes H)d\tau
+\frac{1}{2(2\mu+\lambda)}\int_0^t\|H\|_{L^\infty}^2d\tau\nonumber\\
&\le \frac{1}{2\mu+\lambda}\|(-\Delta)^{-1}\divv(\rho u)\|_{L^\infty}
+\frac{1}{2\mu+\lambda}\|(-\Delta)^{-1}\divv(\rho_0u_0)\|_{L^\infty}\nonumber\\
&\quad+\frac{1}{2\mu+\lambda}\int_0^t\|[u^i, R_{ij}](\rho u^j)\|_{L^\infty}d\tau+
\frac{1}{2\mu+\lambda}\int_0^t\|(-\Delta)^{-1}\divv\divv(H\otimes H)\|_{L^\infty}d\tau\nonumber\\
&\quad+\frac{1}{2(2\mu+\lambda)}\int_0^t\|H\|_{L^\infty}^2d\tau
\triangleq\sum_{i=1}^5Z_i.
\end{align}
By virtue of Gagliardo-Nirenberg inequality, Sobolev's inequality, the $L^q$-estimates of the Riesz transform,
 and H\"older's inequality, we obtain from \eqref{ee1}, \eqref{f5}, and \eqref{f42} that
\begin{align}\label{3.60}
Z_1 &\le \frac{C}{2\mu+\lambda}\|(-\Delta)^{-1}\divv(\rho u) \|_{L^6}^\frac13\|\nabla(-\Delta)^{-1}\divv(\rho u)\|_{L^4}^\frac23\nonumber\\
&\le C\|\rho u\|_{L^2}^\frac13\|\rho u\|_{L^4}^\frac23
\le C\|\rho u\|_{L^2}^\frac12\|\rho u\|_{L^6}^\frac12
\le C(\bar{\rho})\|\sqrt{\rho}u\|_{L^2}^\frac12\|\nabla u\|_{L^2}^\frac12\le C\Big(t_0+\mathbb{E}_0^\frac{4}{5}\Big)^\frac14.
\end{align}
Similarly to $Z_1$, we have
\begin{align}\label{3.61}
Z_2 \le  C\Big(t_0+\mathbb{E}_0^\frac{4}{5}\Big)^\frac14.
\end{align}
For $Z_3$, using the Gagliardo-Nirenberg and Calder{\'o}n-Zygmund inequalities, one deduces that
\begin{align}\label{uuy}
Z_3 & \le \frac{C}{2\mu+\lambda}\int_0^t\|[u^i, R_{ij}](\rho u^j)\|_{L^{12}}^\frac12\|\nabla[u^i, R_{ij}](\rho u^j)\|_{L^4}^\frac12d\tau\nonumber\\
&\le C\int_0^t\big(\|\nabla u\|_{L^3}\|\rho u\|_{L^{12}}\big)^\frac12\big(\|\nabla u\|_{L^6}\|\rho u\|_{L^{12}}\big)^\frac12d\tau\nonumber\\
&\le C(\bar{\rho})\int_0^t\big(\|\nabla u\|_{L^2}+\|\nabla u\|_{L^6}\big)\|u\|_{L^{12}}d\tau\nonumber\\
&\le C(\bar{\rho})\int_0^t\big(\|\nabla u\|_{L^2}+\|\nabla u\|_{L^6}\big)\||u||\nabla u|\|_{L^2}^\frac12d\tau\nonumber\\
&\le C(\bar{\rho})\int_0^t\big(\|\nabla u\|_{L^2}^2+\|\nabla u\|_{L^6}^2+C\|u\|_{L^6}\|\nabla u\|_{L^3}\big)d\tau\nonumber\\
&\le C(\bar{\rho})\int_0^t\Big(\|\nabla u\|_{L^2}^2+\|\nabla u\|_{L^6}^2+\|\nabla u\|_{L^2}^\frac32\|\nabla u\|_{L^6}^\frac12\Big)d\tau\nonumber\\
&\le C(\bar{\rho})\int_0^t\big(\|\nabla u\|_{L^2}^2+\|\nabla u\|_{L^6}^2\big)d\tau\nonumber\\
&\le C(\bar{\rho})\int_0^t\big(\|\nabla u\|_{L^2}^2
+\|\sqrt{\rho}\dot{u}\|_{L^2}^2
+\|H_t\|_{L^2}^2+\|\nabla\theta\|_{L^2}^2\big)d\tau\nonumber\\
&\quad+C(\bar{\rho})\int_0^t\big(\|\nabla H\|_{L^2}^6+\|\nabla u\|_{L^2}^4\|\nabla H\|_{L^2}^2\big)d\tau\nonumber\\
&\le C(\bar{\rho})\left(\int_0^{t_0}+\int_{t_0}^t\right)\big(\|\nabla u\|_{L^2}^2
+\|\sqrt{\rho}\dot{u}\|_{L^2}^2+\|H_t\|_{L^2}^2
+\|\nabla\theta\|_{L^2}^2\big)d\tau
\nonumber\\
&\quad+C\left(\int_0^{t_0}+\int_{t_0}^t\right)\big(\|\nabla H\|_{L^2}^6+\|\nabla u\|_{L^2}^4\|\nabla H\|_{L^2}^2\big)d\tau\nonumber\\
&\le \int_0^t\sigma\big(\|\nabla u\|_{L^2}^2
+\|\sqrt{\rho}\dot{u}\|_{L^2}^2+\|H_t\|_{L^2}^2+\|\nabla\theta\|_{L^2}^2
+\|\nabla H\|_{L^2}^2\big)d\tau+C(M_*)t_0.
\end{align}
Combining \eqref{f5} and \eqref{f42}, one obtains that
\begin{align}\label{sp}
\sup_{0\le t\le T}\big(\|\sqrt{\rho}u\|_{L^2}^2+\|\sqrt{\rho}\theta\|_{L^2}^2+\|H\|_{L^2}^2\big)
+\int_0^T\big(\|\nabla u\|_{L^2}^2+\|\nabla\theta\|_{L^2}^2+\|\nabla H\|_{L^2}^2\big)dt\le C\Big(t_0+\mathbb{E}_0^\frac45\Big).
\end{align}
This together with \eqref{uuy} and \eqref{f42} yields that
\begin{align}\label{f68}
Z_3\le C\Big(t_0+\mathbb{E}_0^\frac45\Big).
\end{align}
For $Z_4$, by H\"older's and Gagliardo-Nirenberg inequalities, \eqref{ee1}, \eqref{w1}, and \eqref{f42}, we have
\begin{align}\label{3.68}
Z_4&\le \frac{1}{2\mu+\lambda}\int_0^t\|(-\Delta)^{-1}\divv\divv(H\otimes H)\|_{L^\infty}d\tau\nonumber\\
&\le C\left(\int_0^{t_0}+\int_{t_0}^t\right)\big(\|\nabla H\|_{L^2}^2+\|\nabla^2H\|_{L^2}^2\big)d\tau\nonumber\\
&\le C(M_*)t_0+\int_0^t\sigma\big(\|\nabla H\|_{L^2}^2+\|\nabla^2H\|_{L^2}^2\big)d\tau\nonumber\\
&\le C\Big(t_0+\mathbb{E}_0^\frac45\Big).
\end{align}
Here we have used the following Gagliardo-Nirenberg inequality
\begin{align}\label{xz}
\|H\|_{L^\infty}\le C\|\nabla H\|_{L^2}^\frac12\|\nabla^2 H\|_{L^2}^\frac12.
\end{align}
Similarly to \eqref{3.68}, we have
\begin{align*}
Z_5\le C\Big(t_0+\mathbb{E}_0^\frac45\Big).
\end{align*}
Substituting \eqref{3.60}, \eqref{3.61}, \eqref{f68}, and \eqref{3.68} into \eqref{xyx}, we find that
\begin{align}
b(t)-b(0)&\le C_5\Big(t_0+\mathbb{E}_0^\frac45\Big)^\frac14\le \ln\frac32,
\end{align}
provided that
$$\mathbb{E}_0\le \varepsilon_4\triangleq\min\left\{\varepsilon_3, \Big(\frac{\ln\frac32}{2^\frac14C_5}\Big)^5\right\}.$$
Integrating \eqref{tb1} with respect to $s$ over $[0, t]$, we get
\begin{align*}
\ln\rho^\delta(x, t)=\ln[\rho_0(X(t; x, 0))+\delta]+\int_0^tg(\tau)d\tau+b(t)-b(0)\le \ln(\bar{\rho}+\delta)+\ln\frac{3}{2}.
\end{align*}
Letting $\delta\rightarrow 0^+$, we arrive at
\begin{align*}
\rho\le \frac{3\bar{\rho}}{2}.
\end{align*}
This finishes the proof of Lemma \ref{l34}.   \hfill $\Box$

{\it Proof of Proposition \ref{p1}.}
1. It follows from \eqref{w1} and \eqref{f42} that
\begin{align}
&\sup_{0\le t\le T}\big(\|\nabla u\|_{L^2}^2
+\|\nabla\theta\|_{L^2}^2+\|\nabla H\|_{L^2}^2+\|\sqrt{\rho}\dot{u}\|_{L^2}^2
+\|H_t\|_{L^2}^2\big)\nonumber\\
&\le \sup_{0\le t\le t_0}\big(\|\nabla u\|_{L^2}^2
+\|\nabla\theta\|_{L^2}^2+\|\nabla H\|_{L^2}^2+\|\sqrt{\rho}\dot{u}\|_{L^2}^2
+\|H_t\|_{L^2}^2\big)\nonumber\\
&\quad+\sup_{0\le t\le T}\big[\sigma\big(\|\nabla u\|_{L^2}^2
+\|\nabla\theta\|_{L^2}^2+\|\nabla H\|_{L^2}^2+\|\sqrt{\rho}\dot{u}\|_{L^2}^2
+\|H_t\|_{L^2}^2\big)\big]\nonumber\\
&\le M_*+C_6\Big(t_0+\mathbb{E}_0^\frac{4}{5}\Big).
\end{align}
Thus, we have
\begin{align}\label{f72}
\sup_{0\le t\le T}\big(\|\nabla u\|_{L^2}^2+\|\nabla\theta\|_{L^2}^2+\|\nabla H\|_{L^2}^2
+\|\sqrt{\rho}\dot{u}\|_{L^2}^2
+\|H_t\|_{L^2}^2\big)\le \frac{3M_*}{2},
\end{align}
provided that
\begin{align*}
t_0+\mathbb{E}_0^\frac{4}{5}\le \frac{M_*}{2C_6}.
\end{align*}

2. From \eqref{ee1}, \eqref{f22}, and \eqref{f50}, we derive that
\begin{align*}
\|\nabla^2\theta\|_{L^2}^2&\le C\big(\|\sqrt{\rho}\dot{\theta}\|_{L^2}^2+\|\theta\nabla u\|_{L^2}^2+\|\nabla u\|_{L^4}^4
+\|\nabla H\|_{L^4}^4\big)\nonumber\\
&\le C\big(\|\sqrt{\rho}\dot{\theta}\|_{L^2}^2+\|\theta\|_{L^\infty}^2\|\nabla u\|_{L^2}^2+\|\nabla u\|_{L^2}\|\nabla u\|_{L^6}^3
+\|\nabla H\|_{L^4}^4\big)\nonumber\\
&\le \frac{1}{2}\|\nabla^2\theta\|_{L^2}^2
+C\|\sqrt{\rho}\dot{\theta}\|_{L^2}^2+C\|\nabla\theta\|_{L^2}^2+C\|\nabla u\|_{L^2}^2
+C\|\sqrt{\rho}\dot{u}\|_{L^2}^2+C\|\nabla H\|_{L^2}\|\nabla^2H\|_{L^2}^3\nonumber\\
&\le  \frac{1}{2}\|\nabla^2\theta\|_{L^2}^2
+C\|\sqrt{\rho}\dot{\theta}\|_{L^2}^2+C\|\nabla\theta\|_{L^2}^2+C\|\nabla u\|_{L^2}^2
+C\|\sqrt{\rho}\dot{u}\|_{L^2}^2+C\|\nabla H\|_{L^2}^2+C\|H_t\|_{L^2}^4,
\end{align*}
which gives that
\begin{align}\label{3.74}
\|\nabla^2\theta\|_{L^2}^2\le C\|\sqrt{\rho}\dot{\theta}\|_{L^2}^2+C\|\nabla\theta\|_{L^2}^2+C\|\nabla u\|_{L^2}^2
+C\|\sqrt{\rho}\dot{u}\|_{L^2}^2+C\|\nabla H\|_{L^2}^2+C\|H_t\|_{L^2}^4.
\end{align}
Thus, using \eqref{w1}, \eqref{f5}, and \eqref{f42}, we have
\begin{align}\label{f74}
&\int_0^T\big(\|\sqrt{\rho}\dot{u}\|_{L^2}^2+\|H_t\|_{L^2}^2
+\|\nabla^2H\|_{L^2}^2+\|\nabla\theta\|_{L^2}^2+\|\nabla^2\theta\|_{L^2}^2\big)dt\nonumber\\
&=\left(\int_0^{t_0}+\int_{t_0}^T\right)\big(\|\sqrt{\rho}\dot{u}\|_{L^2}^2+\|H_t\|_{L^2}^2
+\|\nabla^2H\|_{L^2}^2+\|\nabla\theta\|_{L^2}^2+\|\nabla^2\theta\|_{L^2}^2\big)dt\nonumber\\
&\le M_*t_0+\int_0^T\sigma\big(\|\sqrt{\rho}\dot{u}\|_{L^2}^2+\|H_t\|_{L^2}^2
+\|\nabla^2H\|_{L^2}^2+\|\nabla\theta\|_{L^2}^2+\|\nabla^2\theta\|_{L^2}^2\big)dt\nonumber\\
&\quad+\int_0^T\sigma\big(\|\sqrt{\rho}\dot{\theta}\|_{L^2}^2+\|\nabla H\|_{L^2}^2+\|\nabla u\|_{L^2}^2\big)dt
+C\sup_{0\le t\le T}\|H_t\|_{L^2}^2
\int_0^T\sigma\|H_t\|_{L^2}^2dt\nonumber\\
&\le C_7\Big(t_0+\mathbb{E}_0^\frac45\Big)\le \mathbb{E}_0^\frac23,
\end{align}
provided that
\begin{align*}
t_0+\mathbb{E}_0^\frac{4}{5}\le \frac{\mathbb{E}_0^\frac23}{C_7}.
\end{align*}

Finally, \eqref{ee2} is valid if we select
\begin{align}
\mathbb{E}_0\le \varepsilon_0\quad {\rm and}\quad t_0\le t_1\triangleq\min\left\{\frac{T_*}{2}, \varepsilon_0\right\}
\end{align}
with
\begin{align*}
\varepsilon_0\triangleq\min\left\{\frac{\bar\rho\|\nabla\theta_0\|_{L^2}^6}{54c_v^4\pi^4}, \frac{1}{(2C_1)^\frac32}, \Big(\frac{\delta_1}{8C_3}\Big)^\frac32, \Big(\frac{1}{2C_4}\Big)^\frac32,
\Big(\frac{\log\frac32}{2^\frac14C_5}\Big)^5,
\Big(\frac{M_*}{2C_6}\Big)^\frac54, \Big(\frac{1}{2C_7}\Big)^\frac{15}{2}\right\}.
\end{align*}
The proof of Proposition \ref{p1} is complete. \hfill $\Box$

\begin{lemma}\label{l35}
Under the assumptions in Proposition \ref{p1}, we have
\begin{align}
\sup_{0\le t\le T}\big(\|\sqrt{\rho}\dot{\theta}\|_{L^2}^2+\|\nabla^2\theta\|_{L^2}^2+\|\nabla^2H\|_{L^2}^2\big)
+\int_0^T\|\nabla\dot{\theta}\|_{L^2}^2dt\le C.
\end{align}
\end{lemma}
{\it Proof.}
1. Operating $\partial_t+\divv(u\cdot)$ to $\eqref{a1}_3$ gives rise to
\begin{align}\label{t8}
c_v\rho\big(\dot{\theta}+u\cdot\nabla\dot{\theta}\big)
&=\kappa\Delta\dot{\theta}
+\kappa[\divv u\Delta\theta-
\partial_i(\partial_i\cdot\nabla\theta)-\partial_iu\cdot\nabla\partial_i\theta]
-\rho\dot{\theta}\divv u-\rho\theta\divv\dot{u}\nonumber\\
&\quad+\Big(\lambda(\divv u)^2+\frac{\mu}{2}|\nabla u+(\nabla u)^{tr}|^2\Big)\divv u
+2\lambda\big(\divv\dot{u}-\partial_ku^l\partial_lu^k\big)\divv u\nonumber\\
&\quad+\rho\theta\partial_ku^l\partial_lu^k+\mu\big(\partial_iu^j+\partial_ju^i\big)
\big(\partial_i\dot{u}^j
+\partial_j\dot{u}^i-\partial_iu^k\partial_ku^j-\partial_ju^k\partial_ku^i\big)\nonumber\\
&\quad+\partial_t|\curl H|^2+\divv\big(|\curl H|^2u\big).
\end{align}
Then multiplying \eqref{t8} by $\dot{\theta}$ and integration by parts lead to
\begin{align}\label{t9}
&\frac{c_v}{2}\frac{d}{dt}\int\rho|\dot{\theta}|^2dx+\kappa\int|\nabla\dot{\theta}|^2dx\nonumber\\
&\le C\int|\nabla u|\big(|\nabla^2\theta||\dot{\theta}|+|\nabla\theta||\nabla\dot{\theta}|\big)dx
+C\int|\nabla u|^2|\dot{\theta}|\big(|\nabla u|+\rho\theta\big)dx\nonumber\\
&\quad+C\int\rho|\dot{\theta}|^2|\nabla u|dx+C\int\rho\theta|\nabla\dot{u}||\dot{\theta}|dx
+C\int|\nabla u||\nabla\dot{u}||\dot{\theta}|dx\nonumber\\
&\quad+\int\dot{\theta}\partial_t|\curl H|^2dx+\int\dot{\theta}
\divv\big(|\curl H|^2u\big)dx\triangleq\sum_{i=1}^7K_i.
\end{align}
By virtue of H\"older's, Sobolev's, and Gagliardo-Nirenberg inequalities,
we can bound each term $K_i$ as follows
\begin{align*}
K_1&\le C\|\nabla u\|_{L^3}\big(\|\nabla^2\theta\|_{L^2}\|\dot{\theta}\|_{L^6}
+\|\nabla\dot{\theta}\|_{L^2}\|\nabla\theta\|_{L^6}\big)
\le \frac{\kappa}{14}\|\nabla\dot{\theta}\|_{L^2}^2+C\|\nabla u\|_{L^3}^2\|\nabla^2\theta\|_{L^2}^2,\\
K_2&\le C\|\nabla u\|_{L^3}\|\nabla u\|_{L^4}^2\|\dot{\theta}\|_{L^6}
+C\|\rho\theta\|_{L^3}\|\nabla u\|_{L^4}^2\|\dot{\theta}\|_{L^6}\nonumber\\
&\le \frac{\kappa}{14}\|\nabla\dot{\theta}\|_{L^2}^2
+C\big(\|\rho\theta\|_{L^3}^2+\|\nabla u\|_{L^3}^2\big)\|\nabla u\|_{L^4}^4,\\
K_3&\le C(\bar{\rho})\|\sqrt{\rho}\dot{\theta}\|_{L^2}\|\dot{\theta}\|_{L^6}\|\nabla u\|_{L^3}\le \frac{\kappa}{14}\|\nabla\dot{\theta}\|_{L^2}^2+C\|\sqrt{\rho}\dot{\theta}\|_{L^2}^2\|\nabla u\|_{L^3}^2,\\
K_4&\le C\|\rho\theta\|_{L^3}\|\nabla\dot{u}\|_{L^2}\|\dot{\theta}\|_{L^6}\le \frac{\kappa}{14}\|\nabla\dot{\theta}\|_{L^2}^2
+C\|\rho\theta\|_{L^3}^2\|\nabla\dot{u}\|_{L^2}^2,\\
K_5&\le C\|\nabla u\|_{L^3}\|\nabla\dot{u}\|_{L^2}\|\dot{\theta}\|_{L^6}\le \frac{\kappa}{14}\|\nabla\dot{\theta}\|_{L^2}^2+C\|\nabla u\|_{L^3}^2\|\nabla\dot{u}\|_{L^2}^2,\\
K_6&\le C\int|\nabla H_t||\nabla H||\dot{\theta}|dx\nonumber\\
&\le C\|\nabla H\|_{L^2}^\frac12\|\nabla^2 H\|_{L^2}^\frac12\|\nabla H_t\|_{L^2}\|\dot{\theta}\|_{L^6}\nonumber\\
&\le \frac{\kappa}{12}\|\nabla\dot{\theta}\|_{L^2}^2
+C\|\nabla H\|_{L^2}\|\nabla^2H\|_{L^2}\|\nabla H_t\|_{L^2}^2\nonumber\\
&\le \frac{\kappa}{14}\|\nabla\dot{\theta}\|_{L^2}^2+C\big(\|\nabla H\|_{L^2}^2+\|\nabla^2H\|_{L^2}^2\big)\|\nabla H_t\|_{L^2}^2
\nonumber\\
K_7&\le C\|u\|_{L^6}\|\nabla\dot{\theta}\|_{L^2}\|\nabla H\|_{L^6}^2\le \frac{\kappa}{14}\|\nabla\dot{\theta}\|_{L^2}^2
+C\|\nabla u\|_{L^2}^2\|\nabla^2H\|_{L^2}^4.
\end{align*}
Substituting the above estimates on  $K_i\ (i=1, 2, \ldots, 7)$ into \eqref{t9} yields that
\begin{align}\label{f80}
& c_v\frac{d}{dt}\|\sqrt{\rho}\dot{\theta}\|_{L^2}^2
+\kappa\|\nabla\dot{\theta}\|_{L^2}^2 \notag \\
&\le C\big(\|\sqrt{\rho}\theta\|_{L^3}^2+\|\nabla u\|_{L^3}^2+\|\nabla^2H\|_{L^2}^2+\|\nabla H\|_{L^2}^2\big)
\nonumber\\
&\quad \times
\big(\|\sqrt{\rho}\dot{\theta}\|_{L^2}^2+\|\nabla^2\theta\|_{L^2}^2+\|\nabla H_t\|_{L^2}^2+\|H_t\|_{L^2}^2
+\|\nabla u\|_{L^4}^4+\|\nabla\dot{u}\|_{L^2}^2\big).
\end{align}

2. It follows from \eqref{f5}, \eqref{w1}, \eqref{f50}, and \eqref{f72} that
\begin{align}\label{f81}
&\sup_{0\le t\le T}\big(\|\sqrt{\rho}\theta\|_{L^3}^2+\|\nabla u\|_{L^3}^2+\|\nabla^2H\|_{L^2}^2+\|\nabla H\|_{L^2}^2\big)\nonumber\\
&\le C(\bar{\rho})\sup_{0\le t\le T}\big(\|\sqrt{\rho}\theta\|_{L^2}^2+
\|\nabla \theta\|_{L^2}^2+\|\nabla u\|_{L^2}^2+\|\nabla u\|_{L^6}^2+\|\nabla^2H\|_{L^2}^2+\|\nabla H\|_{L^2}^2\big)\nonumber\\
&\le C(\bar{\rho})\sup_{0\le t\le T}\big(\|\sqrt{\rho}\theta\|_{L^2}^2+
\|\nabla \theta\|_{L^2}^2+\|\nabla u\|_{L^2}^2+\|H_t\|_{L^2}^2+\|\sqrt{\rho}\dot{u}\|_{L^2}^2+\|\nabla H\|_{L^2}^2\big)\nonumber\\
&\le C(M_*, \bar{\rho}).
\end{align}
From \eqref{w1}, \eqref{ee2}, \eqref{f42}, and \eqref{f74}, we get that
\begin{align}\label{f82}
&\int_0^T\big(\|\sqrt{\rho}\dot{\theta}\|_{L^2}^2
+\|\nabla^2\theta\|_{L^2}^2+\|\nabla H_t\|_{L^2}^2+\|H_t\|_{L^2}^2
+\|\nabla u\|_{L^4}^4+\|\nabla\dot{u}\|_{L^2}^2\big)dt\nonumber\\
&\le \left(\int_0^{t_0}+\int_{t_0}^T\right)\big(\|\sqrt{\rho}\dot{\theta}\|_{L^2}^2
+\|\nabla^2\theta\|_{L^2}^2+\|\nabla H_t\|_{L^2}^2+\|H_t\|_{L^2}^2
+\|\nabla u\|_{L^4}^4+\|\nabla\dot{u}\|_{L^2}^2\big)dt\nonumber\\
&\le C(M_*),
\end{align}
where we have used the following fact
\begin{align}\label{fff}
\|\nabla u\|_{L^4}^4&\le \|\nabla u\|_{L^2}\|\nabla u\|_{L^6}^3\nonumber\\
&\le C(M_*)\|\nabla u\|_{L^2}\big(\|\sqrt{\rho}\dot{u}\|_{L^2}+\|H_t\|_{L^2}+\|\nabla\theta\|_{L^2}
+\|\nabla H\|_{L^2}\big)^3\nonumber\\
&\le C\|\nabla u\|_{L^2}^4+C\|\sqrt{\rho}\dot{\theta}\|_{L^2}^4+C\|\nabla H\|_{L^2}^4
+C\|\nabla\theta\|_{L^2}^4+C\|H_t\|_{L^2}^4.
\end{align}
Integrating \eqref{f80} over $[0, T]$, we thus infer from \eqref{f81} and \eqref{f82} that
\begin{align}\label{f83}
\sup_{0\le t\le T}\|\sqrt{\rho}\dot{\theta}\|_{L^2}^2+\int_0^T\|\nabla\dot{\theta}\|_{L^2}^2dt\le C.
\end{align}
This together with \eqref{f83} and \eqref{f72} leads to
\begin{align}
\|\nabla^2\theta\|_{L^2}^2&\le C\|\sqrt{\rho}\dot{\theta}\|_{L^2}^2+C\|\nabla\theta\|_{L^2}^2+C\|\nabla u\|_{L^2}^2
+C\|\sqrt{\rho}\dot{u}\|_{L^2}^2+C\|\nabla H\|_{L^2}^2+C\|H_t\|_{L^2}^4\le C.
\end{align}
Thus, the proof of Lemma \ref{l35} is complete. \hfill $\Box$

As a direct consequence of Proposition \ref{p1} and Lemmas \ref{l31}--\ref{l35}, we arrive at
\begin{corollary}
Under the assumptions in Proposition \ref{p1}, there exists a
positive constant $C$ depending only on $\mu$, $\lambda$, $c_v$, $\kappa$, $\eta$, and the initial data such that the following estimates hold:
\begin{align}\label{tt}
&\sup_{0\le t\le T}\big(\|\nabla u\|_{L^2}^2+\|\nabla\theta\|_{H^1}^2
+\|H\|_{H^2}^2+\|H_t\|_{L^2}^2+\|\sqrt{\rho}\dot{u}\|_{L^2}^2
+\|\sqrt{\rho}\dot{\theta}\|_{L^2}^2\big)\nonumber\\
&\quad+\int_0^T\big(\|\nabla\dot{u}\|_{L^2}^2+\|\nabla\dot{\theta}\|_{L^2}^2+\|\nabla H_t\|_{L^2}^2+\|\nabla u\|_{L^2}^2+\|\nabla\theta\|_{H^1}^2
+\|\nabla H\|_{H^1}^2\big)dt\le C.
\end{align}
\end{corollary}

\subsection{Time-dependent estimates}

\begin{lemma}\label{l36}
Under the same assumptions in Theorem \ref{thm1}, we have
\begin{align}
\sup_{0\le t\le T}\big(\|\nabla\rho\|_{L^2\cap L^q}+\|\nabla^2u\|_{L^2}\big)+\int_0^T\big(\|\nabla^2u\|_{L^q}^2
+\|\nabla^2\theta\|_{L^q}^2
+\|\nabla^3H\|_{L^2}^2\big)dt\le C(T).
\end{align}
\end{lemma}
{\it Proof.}
1. Taking spatial derivative $\nabla$ on the mass equation \eqref{a1}$_1$
leads to
\begin{equation*}
\partial_{t}\nabla\rho+u\cdot\nabla^2\rho
+\nabla u\cdot\nabla\rho+\divv u\nabla\rho+\rho\nabla\divv u=0.
\end{equation*}
For $q\in (3, 6)$, multiplying the above equality by $q|\nabla\rho|^{q-2}\nabla\rho$ gives that
\begin{align*}
(|\nabla\rho|^q)_t+\divv(|\nabla\rho|^qu)+(q-1)|\nabla\rho|^q\divv u+q|\nabla\rho|^{q-2}(\nabla\rho)^{tr}\nabla u(\nabla\rho)+q\rho|\nabla\rho|^{q-2}\nabla\rho\cdot\nabla\divv u=0.
\end{align*}
Thus, integration by parts over $\mathbb{R}^3$ yields that
\begin{align}\label{0t86}
\frac{d}{dt}\|\nabla\rho\|_{L^q}
\le C\big(\|\nabla u\|_{L^\infty}+1\big)\|\nabla\rho\|_{L^q}+\|\nabla^2u\|_{L^q}.
\end{align}
Applying the standard $L^q$-estimate to \eqref{a1}$_2$ leads to
\begin{align}\label{t87}
\|\nabla^2u\|_{L^q}
\le C\big(\|\rho\dot{u}\|_{L^q}+\|\nabla (\rho\theta)\|_{L^q}+\||H||\nabla H|\|_{L^q}\big)
\le C\big(1+\|\nabla\dot{u}\|_{L^2}+\|\nabla\rho\|_{L^q}\big),
\end{align}
due to \eqref{tt}. So we obtain from \eqref{0t86} and \eqref{t87} that
\begin{align}\label{t86}
\frac{d}{dt}\|\nabla\rho\|_{L^q}\le C\big(\|\nabla u\|_{L^\infty}+1\big)\|\nabla\rho\|_{L^q}+C\big(1+\|\nabla\dot{u}\|_{L^2}
+\|\nabla\rho\|_{L^q}\big).
\end{align}

2.  It follows from the Gagliardo-Nirenberg inequality and \eqref{w5} that
\begin{align}\label{t88}
\|\divv u\|_{L^\infty}+\|\curl u\|_{L^\infty}
&\le C\big(\|F\|_{L^\infty}+\|\rho\theta\|_{L^\infty}
+\|w\|_{L^\infty}\big)\nonumber\\
&\le C\|F\|_{L^6}^\frac12\|\nabla F\|_{L^6}^\frac12+C\|w\|_{L^6}^\frac12\|\nabla w\|_{L^6}^\frac12+C\|\nabla\theta\|_{L^2}^\frac12\|\nabla^2\theta\|_{L^2}^\frac12\nonumber\\
&\le C\big(1+\|\nabla F\|_{L^2}+\|\nabla F\|_{L^6}+\|\nabla w\|_{L^2}+\|\nabla w\|_{L^6}\big)\nonumber\\
&\le C\big(1+\|\rho\dot{u}\|_{L^2}
+\||H||\nabla H|\|_{L^2}+\|\rho\dot{u}\|_{L^6}
+\||H||\nabla H|\|_{L^6}\big)\nonumber\\
&\le C(\bar{\rho})\big(1+\|\sqrt{\rho}\dot{u}\|_{L^2}+\|\nabla\dot{u}\|_{L^2}
+\|\nabla^2H\|_{L^2}^2\big)\nonumber\\
&\le C(\bar{\rho})\big(1+\|\nabla\dot{u}\|_{L^2}\big),
\end{align}
which together with Lemma \ref{l24}, \eqref{t87}, and \eqref{t88} indicates that
\begin{align}\label{t89}
\|\nabla u\|_{L^\infty}
&\le C\big(\|\divv u\|_{L^\infty}+\|\curl u\|_{L^\infty}\big)
\ln\big(e+\|\nabla^2u\|_{L^q}\big)
+C\|\nabla u\|_{L^2}+C\nonumber\\
&\le C\big(1+\|\nabla\dot{u}\|_{L^2}\big)\ln\big(
e+\|\nabla\dot{u}\|_{L^2}+\|\nabla\rho\|_{L^q}\big)+C.
\end{align}
Substituting \eqref{t89} into \eqref{t86} yields that
\begin{align}\label{t90}
f'(t)\leq Cg(t)f(t)\ln f(t),
\end{align}
where
\begin{align*}
f(t)\triangleq e+\|\nabla\rho\|_{L^q},\ g(t)\triangleq 1+\|\nabla\dot{u}\|_{L^2}.
\end{align*}
It thus follows from \eqref{t90}, Gronwall's inequality, and \eqref{tt} that
\begin{align}\label{t91}
\sup_{0\le t\le T}\|\nabla\rho\|_{L^q}\le C(T).
\end{align}
It should be noted that \eqref{0t86} also holds true for $q=2$, hence we then derive from Gronwall's inequality, \eqref{t89}, \eqref{t91}, \eqref{t87},
and \eqref{tt} that
\begin{align}\label{zt91}
\sup_{0\le t\le T}\|\nabla\rho\|_{L^2}\le C(T).
\end{align}
This combined with \eqref{t87}, \eqref{f72}, \eqref{f81}, and \eqref{f82} implies that
\begin{align}\label{t92}
\sup_{0\le t\le T}\|\nabla^2u\|_{L^2}&\le C\sup_{0\le t\le T}\big(\|\sqrt{\rho}\dot{u}\|_{L^2}
+\|\nabla(\rho\theta)\|_{L^2}+\||H||\nabla H|\|_{L^2}\big)\nonumber\\
&\le C+C\sup_{0\le t\le T}\big(\|\nabla\rho\|_{L^2}+\|\nabla\theta\|_{L^2}+\|\nabla^2H\|_{L^2}\big)\nonumber\\
&\le C(T),
\end{align}
and
\begin{align}
\int_0^T\|\nabla^2u\|_{L^q}^2dt\le C\int_0^T\big(1+\|\nabla\dot{u}\|_{L^2}^2+\|\nabla\rho\|_{L^q}^2\big)dt\le C(T).
\end{align}

3. Taking the operator $\nabla$ to $\eqref{a1}_4$, we get
\begin{align}\label{3.97}
-\eta\nabla\Delta H=\nabla\big(H\cdot\nabla u-u\cdot\nabla H-H\divv u-H_t\big).
\end{align}
Using the $L^2$-estimate of elliptic equations, we derive from \eqref{3.97} that
\begin{align}
\|\nabla^3H\|_{L^2}&\le C\|\nabla(H\cdot\nabla u-u\cdot\nabla H-H\divv u-H_t)\|_{L^2}\nonumber\\
&\le C\|\nabla H_t\|_{L^2}+C\||\nabla^2u||H|\|_{L^2}+C\||\nabla u||\nabla H|\|_{L^2}
+C\||u||\nabla^2H|\|_{L^2}\nonumber\\
&\le C\|\nabla H_t\|_{L^2}+C\|H\|_{L^\infty}\|\nabla^2u\|_{L^2}+C\|\nabla u\|_{L^3}\|\nabla H\|_{L^6}
+C\|u\|_{L^6}\|\nabla^2H\|_{L^3}\nonumber\\
&\le C\|\nabla H\|_{L^2}^\frac{1}{2}\|\nabla^2H\|_{L^2}^\frac{1}{2}\|\nabla^2u\|_{L^2}
+C\|\nabla u\|_{L^2}^\frac{1}{2}\|\nabla^2u\|_{L^2}^\frac{1}{2}\|\nabla^2H\|_{L^2}\nonumber\\
&\quad+C(\|\nabla u\|_{L^2}^2+\|\nabla H\|_{L^2}^2)\|\nabla^2H\|_{L^2}
+C\|\nabla H_t\|_{L^2}+\frac{1}{2}\|\nabla^3H\|_{L^2},
\end{align}
which yields
\begin{align}
\|\nabla^3H\|_{L^2}^2&\le C\|\nabla H\|_{L^2}\|\nabla^2H\|_{L^2}\|\nabla^2u\|_{L^2}^2
+C\|\nabla u\|_{L^2}\|\nabla^2 u\|_{L^2}\|\nabla^2H\|_{L^2}^2\nonumber\\
&\quad+C(\|\nabla u\|_{L^2}^4+\|\nabla H\|_{L^2}^4)\|\nabla^2H\|_{L^2}^2
+C\|\nabla H_t\|_{L^2}^2,
\end{align}
from which, \eqref{t92}, and \eqref{f72}, one obtains that
\begin{align}\label{t97}
\int_0^T\|\nabla^3H\|_{L^2}^2dt&\le C\int_0^T\big(\|\nabla^2H\|_{L^2}^2+\|\nabla H\|_{L^2}^2+\|\nabla H_t\|_{L^2}^2\big)dt
\le C(T).
\end{align}
This together with \eqref{tt}, \eqref{f72}, \eqref{t92}, and \eqref{t97} leads to
\begin{align}
\int_0^T\|\nabla^2\theta\|_{L^q}^2dt&\le C\int_0^T\big(\|\rho\dot{\theta}\|_{L^q}^2
+\|\theta\nabla u\|_{L^q}^2+\|\nabla u\|_{L^{2q}}^4
+\|\nabla H\|_{L^{2q}}^4\big)dt\nonumber\\
&\le C\int_0^T\big(1+\|\nabla\dot{\theta}\|_{L^2}^2
+\|\theta\|_{L^\infty}^2\|\nabla u\|_{H^1}^2
+\|\nabla u\|_{L^3}^2\|\nabla^2u\|_{L^q}^2
+\|\nabla H\|_{L^3}^2\|\nabla^2H\|_{L^q}^2\big)dt\nonumber\\
&\le C\int_0^T\big(1+\|\nabla\dot{\theta}\|_{L^2}^2
+\|\nabla^2u\|_{L^q}^2+\|\nabla^3H\|_{L^2}^2\big)dt\nonumber\\
&\le C(T).
\end{align}
The proof of Lemma \ref{l36} is finished. \hfill $\Box$

\section{Proof of Theorem \ref{thm1}}
In terms of Lemma \ref{l21}, the Cauchy problem \eqref{a1}--\eqref{a3} has a unique solution $(\rho, u, \theta, H)$ in $\mathbb{R}^3\times(0, T_*]$
for some small $T_*>0$. With the \textit{a priori} estimates obtained in Section 3, we can extend the $(\rho, u, \theta, H)$ to all time. In fact, by the
continuity of density and Lemma \ref{l21}, we see that \eqref{ee1} and \eqref{ee2} are valid in $[0, T_1]$ with
$T_1\triangleq\min\left\{T_*,\frac{2\mathbb{E}_0^\frac23}{M_*}\right\}$.

Set
\begin{align}
T^*=\sup\{T~|~\eqref{ee1}~is~valid\}.
\end{align}
Then $T^*\ge T_1$. We next claim that
\begin{align}
T^*=\infty.
\end{align}
Suppose, by contradiction,
that $T^*<\infty$. Note that, all
the {\it a priori} estimates obtained in Section 3 are uniformly bounded for any $t<T^*$. Hence,
we define
\begin{align}
(\rho, u, \theta, H)(x, T^*)=\lim_{t\rightarrow T^*}(\rho, u, \theta, H)(x, t).
\end{align}
Furthermore, standard arguments yield that $\rho\dot{u}$, $\rho\dot{\theta}\in C([0, T^*); L^2)$, which implies that
\begin{align}
(\rho\dot{u}, \rho\dot{\theta})(x, T^*)=\lim_{t\rightarrow T^*}(\rho\dot{u}, \rho\dot{\theta})(x, t)\in L^2.
\end{align}
Hence,
\begin{align*}
\big[-\mu\Delta u-(\mu+\lambda)\nabla\divv u+\nabla(R\rho\theta)+\curl H\times H\big]|_{t=T^*}&=\sqrt{\rho}(x, T^*)g_1(x),\\
\big[\kappa\Delta \theta+\lambda(\divv u)^2+\frac{\mu}{2}|\nabla u+(\nabla u)^{tr}|^2+\eta|\curl H|^2\big]|_{t=T^*}&=\sqrt{\rho}(x, T^*)g_2(x),
\end{align*}
with
\begin{align*}
g_1(x)=\left\{
\begin{array}{ll}
\displaystyle
\rho^{-\frac{1}{2}}(x, T^*)(\rho\dot{u})(x, T^*), &{\rm for}~ x\in\{x|\rho(x, T^*)>0\},\\[3pt]
0,  &{\rm for}~ x\in\{x|\rho(x, T^*)=0\},
\end{array}
\right.
\end{align*}
and
\begin{align*}
g_2(x)=\left\{
\begin{array}{ll}
\displaystyle
\rho^{-\frac{1}{2}}(x, T^*)(c_v\rho\dot{\theta}+R\rho\theta\divv u)(x, T^*), &{\rm for}~ x\in\{x|\rho(x, T^*)>0\},\\[3pt]
0,  &{\rm for}~ x\in\{x|\rho(x, T^*)=0\},
\end{array}
\right.
\end{align*}
satisfying $(g_1, g_2)\in L^2$ due to \eqref{tt}. Thus, $(\rho, u, \theta, H)$  satisfies compatibility conditions \eqref{qqw1} and \eqref{qqw}. Therefore, we can apply Lemma \ref{l21} to extend the local strong solution
beyond $T^*$ by taking $(\rho, u, \theta, H)(x, T^*)$ as the initial data.
 This contradicts the assumption of $T^*$. This finishes the proof of Theorem
\ref{thm1}.
\hfill $\Box$

\section{Proof of Theorem \ref{thm2}}\label{sec5}

We divide the proof of Theorem \ref{thm2} into several steps.

\textbf{Step 1. Estimate for $\|\nabla H(t)\|_{L^2}$}.
It follows from \eqref{f20} that
\begin{align*}
\frac{d}{dt}\|\nabla H\|_{L^2}^2+\|H_t\|_{L^2}^2+\|\nabla^2 H\|_{L^2}^2
\le C\|\nabla u\|_{L^2}^4\|\nabla H\|_{L^2}^2,
\end{align*}
which implies that
\begin{align*}
\frac{d}{dt}\big(t\|\nabla H\|_{L^2}^2\big)+t\|H_t\|_{L^2}^2
+t\|\nabla^2 H\|_{L^2}^2
&\le \|\nabla H\|_{L^2}^2+C\|\nabla u\|_{L^2}^2\|\nabla u\|_{L^2}^2\big(t\|\nabla H\|_{L^2}^2\big).
\end{align*}
This along with Gronwall's inequality, \eqref{sp}, and \eqref{tt} yields that
\begin{align}\label{dd1}
\sup_{0\le t\le T}\big(t\|\nabla H\|_{L^2}^2\big)
+\int_0^Tt\left(\|H_t\|_{L^2}^2+\|\nabla^2H\|_{L^2}^2\right)dt\le C.
\end{align}

\textbf{Step 2. Estimate for $\|\sqrt{\rho}u(t)\|_{L^2}$ and $\|\sqrt{\rho}\theta(t)\|_{L^2}$}. Multiplying $\eqref{a1}_2$ by $u$ and integrating by parts, we obtain from H\"older's inequality, Sobolev's inequality, \eqref{lz3.1}, \eqref{ee1}, and \eqref{tt} that
\begin{align*}
&\frac12\frac{d}{dt}\|\sqrt{\rho}u\|_{L^2}^2
+\mu\|\nabla u\|_{L^2}^2+(\mu+\lambda)\|\divv u\|_{L^2}^2\nonumber\\
&= \int R\rho\theta\divv udx +\int \curl H\times H \cdot udx\nonumber\\
&\le (\mu+\lambda)\|\divv u\|_{L^2}^2
+C\|\rho\|_{L^\infty}\|\sqrt{\rho}\theta\|_{L^2}^2
+\|\nabla H\|_{L^2}\|H\|_{L^3}\|u\|_{L^6}\nonumber\\
&\le (\mu+\lambda)\|\divv u\|_{L^2}^2
+C(\bar{\rho})\|\rho\|_{L^\frac{3}{2}}\|\theta\|_{L^6}^2
+C\|\nabla H\|_{L^2}^{\frac32}\|H\|_{L^2}^{\frac12}\|\nabla u\|_{L^2}\nonumber\\
&\le (\mu+\lambda)\|\divv u\|_{L^2}^2+c_5\|\nabla\theta\|_{L^2}^2
+\frac{\mu}{2}\|\nabla u\|_{L^2}^2+C\|\nabla H\|_{L^2}^3.
\end{align*}
This gives that
\begin{align}\label{x1}
&\frac{d}{dt}\|\sqrt{\rho}u\|_{L^2}^2
+\mu\|\nabla u\|_{L^2}^2
\leq 2c_5\|\nabla\theta\|_{L^2}^2+C\|\nabla H\|_{L^2}^3.
\end{align}
Multiplying $\eqref{a1}_3$ by $\theta$ and integration by parts, we infer from Sobolev's inequality, Gagliardo-Nirenberg inequality, \eqref{tt}, and \eqref{f42} that
\begin{align*}
&\frac{c_v}{2}\frac{d}{dt}\int\rho\theta^2dx
+\kappa\int|\nabla\theta|^2dx\nonumber\\
&=\int\theta\Big(\lambda(\divv u)^2+\frac{\mu}{2}|\nabla u+(\nabla u)^{tr}|^2-R\rho\theta\divv u\Big)dx+\eta\int(\curl H)\cdot(\curl H)\theta dx\nonumber\\
&\le C\int\theta|\nabla u|^2dx+C\int\rho\theta^2|\nabla u|dx
+C\int|H||\nabla^2H|\theta dx
+C\int|\nabla\theta||H||\nabla H|dx\nonumber\\
&\le C\|\theta\|_{L^\infty}\|\nabla u\|_{L^2}^2
+C\|\theta\|_{L^\infty}\|\nabla u\|_{L^2}\|\theta\|_{L^6}\|\rho\|_{L^3}
+C\|\nabla\theta\|_{L^2}\|\nabla^2H\|_{L^2}\|H\|_{L^3}\nonumber\\
&\le C\|\theta\|_{L^\infty}\big(\|\nabla u\|_{L^2}^2
+\|\nabla\theta\|_{L^2}^2\big)
+C\|\nabla\theta\|_{L^2}\|\nabla^2H\|_{L^2}\|H\|_{L^2}^{\frac12}
\|\nabla H\|_{L^2}^{\frac12}\nonumber\\
&\le C\|\nabla\theta\|_{L^2}^\frac12\|\nabla^2\theta\|_{L^2}^\frac12\big(\|\nabla u\|_{L^2}^2+\|\nabla\theta\|_{L^2}^2\big)
+\frac{\kappa}{2}\|\nabla\theta\|_{L^2}^2
+C\|\nabla H\|_{L^2}\|\nabla^2H\|_{L^2}^2\nonumber\\
&\le C\Big(t_0+\mathbb{E}_0^\frac45\Big)^\frac12\big(\|\nabla u\|_{L^2}^2+\|\nabla\theta\|_{L^2}^2\big)
+\frac{\kappa}{2}\|\nabla\theta\|_{L^2}^2
+C\|\nabla H\|_{L^2}\|\nabla^2H\|_{L^2}^2,
\end{align*}
where we have used the following fact
\begin{align*}
\|\nabla\theta\|_{L^2}\le C\sup_{t_0\le t\le T}\big(\sigma\|\nabla\theta\|_{L^2}^2\big)^\frac12
\le C\Big(t_0+\mathbb{E}_0^\frac45\Big)^\frac12\ \ \text{for}\ \ t\ge t_0.
\end{align*}
Hence, we have
\begin{align}\label{x2}
c_v\frac{d}{dt}\|\sqrt{\rho}\theta\|_{L^2}^2+\kappa\|\nabla\theta\|_{L^2}^2
\le C\|\nabla H\|_{L^2}\|\nabla^2H\|_{L^2}^2
+C_9\Big(t_0+\mathbb{E}_0^\frac45\Big)^\frac12\big(\|\nabla u\|_{L^2}^2+\|\nabla\theta\|_{L^2}^2\big).
\end{align}
Adding \eqref{x2} multiplied by $\Big(\frac{2c_5}{\kappa}+1\Big)$ to $\eqref{x1}$ and choosing
$$\mathbb{E}_0\le \varepsilon_4\triangleq\min\left\{\varepsilon_0,
\Big(\frac{\mu}{2^\frac52C_9}\Big)^{10}, \Big(\frac{\kappa}{2^\frac52C_9}\Big)^{10}\right\},$$
we obtain that
\begin{align}\label{x4} \frac{d}{dt}\big(\|\sqrt{\rho}u\|_{L^2}^2+\|\sqrt{\rho}\theta\|_{L^2}^2\big)
+\|\nabla u\|_{L^2}^2+\|\nabla \theta\|_{L^2}^2
\le C\|\nabla H\|_{L^2}\|\nabla H\|_{H^1}^2.
\end{align}
Multiplying \eqref{x4} by $t^\frac12$ and integrating over $[0, T]$, we obtain from \eqref{tt} and \eqref{dd1} that
\begin{align}\label{5.5}
&\sup_{0\le t\le T}\big[t^\frac12\big(\|\sqrt{\rho}u\|_{L^2}^2+\|\sqrt{\rho}\theta\|_{L^2}^2\big)\big]
+\int_0^Tt^\frac12\big(\|\nabla u\|_{L^2}^2+\|\nabla\theta\|_{L^2}^2\big)dt\nonumber\\
&\le \sup_{0\le t\le T}\big(\|\sqrt{\rho}u\|_{L^2}^2+\|\sqrt{\rho}\theta\|_{L^2}^2\big)\int_0^{t_0}t^{-\frac12}dt
+C\int_{t_0}^T\big(\|\nabla u\|_{L^2}^2+\|\nabla\theta\|_{L^2}^2\big)dt\nonumber\\
&\quad+\sup_{0\le t\le T}\big(t\|\nabla H\|_{L^2}^2\big)^\frac12\int_0^T\|\nabla H\|_{H^1}^2dt
\le C.
\end{align}

\textbf{Step 3. Estimate for $\|\nabla u(t)\|_{L^2}$},
\textbf{$\|\nabla\theta(t)\|_{L^2}$}, \textbf{and $\|\sqrt{\rho}\dot{u}(t)\|_{L^2}$}. Multiplying \eqref{f12} by $t^\frac12$ and integrating the resulting inequality over $[0, T]$, we have
\begin{align}\label{k6}
& t^\frac12\int\left[\frac{\mu}{2}|\nabla u|^2+\frac{\kappa\delta}{2}|\nabla\theta|^2+\eta|\nabla H|^2
+\Big(\frac{2c_2}{\mu}+1\Big)\delta^\frac{3}{2}\rho|\dot{u}|^2
+2\delta^\frac32|H_t|^2\right]dx\nonumber\\
&\quad+\int_0^T\int t^\frac12\Big(\rho|\dot{u}|^2+\frac12|H_t|^2+\frac{\eta^2}{4}|\nabla^2H|^2
+\frac{c_v\delta}{4}\rho|\dot{\theta}|^2
+\mu\delta^\frac32|\nabla \dot{u}|^2+\frac{\eta\delta^\frac32}{2}|\nabla H_t|^2\Big)dxdt\nonumber\\
&\le \int_0^T\int t^{-\frac12}\Big(\frac{\mu}{2}|\nabla u|^2+\frac{\kappa\delta}{2}|\nabla\theta|^2+
\eta|\nabla H|^2+\Big(\frac{2c_2}{\mu}+1\Big)\delta^\frac32\rho|\dot{u}|^2
+2\delta^\frac32|H_t|^2\Big)dxdt\nonumber\\
&\quad+t^\frac12\int\Big(R\rho\theta\divv u+\frac12|H|^2\divv u-H\cdot\nabla u\cdot H
+\lambda\delta\theta(\divv u)^2+\frac{\mu\delta}{2}\theta|\nabla u+(\nabla u)^{tr}|^2\Big)dx\nonumber\\
&\quad+\int_0^T\int t^{-\frac12}\Big(R\rho\theta\divv u+\frac12|H|^2
\divv u-H\cdot\nabla u\cdot H
+\lambda\delta\theta(\divv u)^2+\frac{\mu\delta}{2}\theta|\nabla u+(\nabla u)^{tr}|^2\Big)dxdt\nonumber\\
&\quad +\delta t^{\frac12}\int\theta|\curl H|^2dx
+\delta\int_0^T\int t^{-\frac12}\theta|\nabla H|^2dxdt
+C\delta\int_0^Tt^\frac12\|\theta\|_{L^\infty}^2\|\nabla H\|_{L^2}^2dt \notag \\
& \quad +\delta^{-\frac32}\int_0^Tt^\frac12\big(\|\theta\nabla u\|_{L^2}^2+\|\nabla u\|_{L^2}^2+\|\nabla u\|_{L^2}^4\|\nabla\theta\|_{L^2}^2\big)dt
\nonumber\\
&\quad+C\int_0^Tt^\frac12\big(\delta^\frac32\|\nabla u\|_{L^4}^4+\delta^\frac52\|\nabla^2H\|_{L^2}^4
+\|\nabla H\|_{L^2}^4+\|\nabla u\|_{L^2}^4\|\nabla H\|_{L^2}^2
+\|H\|_{L^2}\|\nabla H\|_{L^2}^3\big)dt.
\end{align}
It follows from \eqref{w1}, \eqref{f5}, and \eqref{tt} that
\begin{align}\label{k7}
&\int_0^T\int t^{-\frac12}\Big(\frac{\mu}{2}|\nabla u|^2+\frac{\kappa\delta}{2}|\nabla\theta|^2+
\eta|\nabla H|^2+\Big(\frac{2c_2}{\mu}+1\Big)\delta^\frac32\rho|\dot{u}|^2
+2\delta^\frac32|H_t|^2\Big)dxdt\nonumber\\
&\le C\sup_{0\le t\le t_0}\big(\|\nabla u\|_{L^2}^2+\|\nabla\theta\|_{L^2}^2+\|\nabla H\|_{L^2}^2
+\|\sqrt{\rho}\dot{u}\|_{L^2}^2+\|H_t\|_{L^2}^2\big)\int_{0}^{t_0}t^{-\frac12}dt\nonumber\\
&\quad+C\int_{t_0}^T\big(\|\nabla u\|_{L^2}^2+\|\nabla\theta\|_{L^2}^2+\|\nabla H\|_{L^2}^2
+\|\sqrt{\rho}\dot{u}\|_{L^2}^2+\|H_t\|_{L^2}^2\big)dt
\le C.
\end{align}
Moreover, we deduce from \eqref{3.49}, \eqref{5.5}, and \eqref{tt} that
\begin{align}
& t^\frac12\int\Big(R\rho\theta{\rm div}u+\frac12|H|^2{\rm div}u-H\cdot\nabla u\cdot H
+\lambda\delta\theta({\rm div}u)^2+\frac{\mu\delta}{2}\theta|\nabla u+(\nabla u)^{tr}|^2+\delta\theta|\nabla H|^2\Big)dx \notag \\
&\leq Ct^\frac12\int\big(\rho\theta|{\rm div}u|
+|\nabla u||H|^2+\delta\theta|\nabla H|^2+\delta\theta|\nabla u|^2\big)dx \nonumber\\
&\le C(\bar\rho)t^\frac12\|\sqrt{\rho}\theta\|_{L^2}\|\nabla u\|_{L^2}
+Ct^\frac12\|\nabla u\|_{L^2}\|H\|_{L^4}^2+C\delta t^\frac12\|\theta\|_{L^\infty}
\big(\|\nabla H\|_{L^2}^2+\|\nabla u\|_{L^2}^2\big)\nonumber\\
&\le C\delta^{-1}t^\frac12\|\sqrt{\rho}\theta\|_{L^2}^2
+C\delta t^\frac12\|\nabla u\|_{L^2}^2
+Ct^\frac12\|\nabla u\|_{L^2}\|H\|_{L^2}^{\frac12}\|\nabla H\|_{L^2}^{\frac32}
\notag \\
&\le C+C\delta\sup_{0\le t\le T}\big[t^\frac12\big(\|\nabla H\|_{L^2}^2
+\|\nabla u\|_{L^2}^2\big)\big].
  \end{align}
We obtain from \eqref{w1}, H\"older's inequality, \eqref{lz3.1}, \eqref{3.49}, Gagliardo-Nirenberg inequality,
and \eqref{tt} that
\begin{align}
& \int_0^T\int t^{-\frac12}\Big(R\rho\theta{\rm div}u+\frac12|H|^2{\rm div}u-H\cdot\nabla u\cdot H
+\lambda\delta\theta({\rm div}u)^2+\frac{\mu\delta}{2}\theta|\nabla u+(\nabla u)^{tr}|^2+\delta\theta|{\rm curl} H|^2\Big)dxdt \notag \\
& \leq C\int_0^T\int t^{-\frac12}\big(\rho\theta|{\rm div}u|+\theta|\nabla u|^2+|\nabla u||H|^2+\theta|\nabla H|^2\big)dxdt\nonumber\\
&= C\left(\int_0^{t_0}+\int_{t_0}^T\right)\int t^{-\frac12}\big(\rho\theta|{\rm div}u|+\theta|\nabla u|^2+|\nabla u||H|^2+\theta|\nabla H|^2\big)dxdt\nonumber\\
&\le C\sup_{0\le t\le t_0}\big(\|\sqrt{\rho}\theta\|_{L^2}^2+\|\nabla\theta\|_{L^2}^2
+\|\nabla^2\theta\|_{L^2}^2
+\|\nabla u\|_{L^2}^2+\|\nabla H\|_{L^2}^2\big)\int_0^{t_0}t^{-\frac12}dt\nonumber\\
&\quad+C\int_{t_0}^T\big(\|\sqrt{\rho}\theta\|_{L^2}\|\nabla u\|_{L^2}
+\|\theta\|_{L^\infty}\|\nabla u\|_{L^2}^2
+\|\nabla u\|_{L^2}\|H\|_{L^4}^2+\|\theta\|_{L^\infty}\|\nabla H\|_{L^2}^2\big)dt\nonumber\\
&\le C(t_0,\|\rho_0\|_{L^1},\bar\rho)\sup_{0\le t\le t_0}\big(\|\nabla\theta\|_{H^1}^2
+\|\nabla u\|_{L^2}^2+\|\nabla H\|_{L^2}^2\big)
+C\int_{t_0}^T\big(\|\nabla u\|_{L^2}^2+\|\nabla\theta\|_{L^2}^2\big)dt
\nonumber\\
&\quad +C\sup_{t_0\le t\le T}\Big(\|\nabla\theta\|_{L^2}^{\frac12}
\|\nabla^2\theta\|_{L^2}^{\frac12}\Big)
\int_{t_0}^T
\big(\|\nabla u\|_{L^2}^2+\|\nabla H\|_{L^2}^2\big)dt
\nonumber\\
& \quad +C\sup_{t_0\le t\le T}\|\nabla u\|_{L^2}^2
\int_{t_0}^T\|\nabla u\|_{L^2}^2dt
+C\int_{t_0}^T\|\nabla H\|_{L^2}^2dt\le C.
\end{align}
By Gagliardo-Nirenberg inequality, \eqref{tt}, and \eqref{dd1}, we get that
\begin{align}
\int_0^Tt^\frac12\|\theta\|_{L^\infty}^2\|\nabla H\|_{L^2}^2dt
& \leq C\int_0^Tt^\frac12\|\nabla\theta\|_{L^2}\|\nabla^2\theta\|_{L^2}\|\nabla H\|_{L^2}^2dt \notag \\
& \le C\sup_{0\le t\le T}\big(t\|\nabla H\|_{L^2}^2\big)^\frac12\int_{0}^T\big(\|\nabla\theta\|_{L^2}^2
+\|\nabla H\|_{L^2}^2\big)dt
\le C.
\end{align}
We infer from Gagliardo-Nirenberg inequality, \eqref{tt}, and \eqref{5.5} that
\begin{align}
&\int_0^Tt^\frac12\big(\|\theta\nabla u\|_{L^2}^2+\|\nabla u\|_{L^2}^2+\|\nabla u\|_{L^2}^4\|\nabla\theta\|_{L^2}^2\big)dx\nonumber\\
&\le C\sup_{0\le t\le T}\Big(\|\nabla \theta\|_{L^2}^{\frac12}
\|\nabla^2\theta\|_{L^2}^{\frac12}+\|\nabla u\|_{L^2}^4+1\Big)
\int_0^Tt^\frac12\big(\|\nabla u\|_{L^2}^2+\|\nabla\theta\|_{L^2}^2\big)dt
\le C.
\end{align}
From \eqref{w7}, \eqref{xz}, \eqref{tt}, and Young's inequality, we have
\begin{align}\label{5.17}
\|\nabla u\|_{L^4}^4&\le \|\nabla u\|_{L^2}\|\nabla u\|_{L^6}^3\nonumber\\
&\le C\|\nabla u\|_{L^2}\big(\|\sqrt{\rho}\dot{u}\|_{L^2}+\|\nabla\theta\|_{L^2}
+\|H\|_{L^\infty}\|\nabla H\|_{L^2}\big)^3\nonumber\\
&\le C\|\nabla u\|_{L^2}^4+C\|\sqrt{\rho}\dot{u}\|_{L^2}^4
+C\|\nabla\theta\|_{L^2}^4+C\|\nabla H\|_{L^2}^2\|\nabla^2 H\|_{L^2}^2.
\end{align}
Hence, we obtain from \eqref{tt}, \eqref{5.5}, and Proposition \ref{p1} that
\begin{align}
\int_0^Tt^\frac12\|\nabla u\|_{L^4}^4dt
&\le C\int_0^Tt^\frac12\big(
\|\nabla u\|_{L^2}^4+\|\sqrt{\rho}\dot{u}\|_{L^2}^4
+\|\nabla\theta\|_{L^2}^4+\|\nabla H\|_{L^2}^2\|\nabla^2 H\|_{L^2}^2\big)dt\nonumber\\
&\le C\sup_{0\le t\le T}
\big(\|\nabla u\|_{L^2}^2+\|\nabla \theta\|_{L^2}^2\big)
\int_0^Tt^\frac12\big(\|\nabla u\|_{L^2}^2
+\|\nabla\theta\|_{L^2}^2\big)dt
\nonumber\\
&\quad +C\sup_{0\le t\le T}
\big[t^\frac12\big(\|\sqrt{\rho}\dot{u}\|_{L^2}^2+\|\nabla H\|_{L^2}^2\big)\big]
\int_0^T\big(\|\sqrt{\rho}\dot{u}\|_{L^2}^2+\|\nabla^2 H\|_{L^2}^2\big)dt \notag \\
& \leq C+C_{10}\mathbb{E}_0^\frac23\sup_{0\le t\le T}\big[t^\frac12(\|\sqrt{\rho}\dot{u}\|_{L^2}^2+\|\nabla H\|_{L^2}^2)\big].
\end{align}
From \eqref{f49}, \eqref{tt}, and Proposition \ref{p1}, we have
\begin{align}
\int_0^Tt^\frac12\|\nabla^2H\|_{L^2}^4ds
&\le C\sup_{0\le t\le T}\big[t^\frac12(\|H_t\|_{L^2}^2
+\|\nabla H\|_{L^2}^2)\big]\int_0^T\|\nabla^2H\|_{L^2}^2dt\nonumber\\
&\le C_{11}\mathbb{E}_0^\frac23\sup_{0\le t\le T}\big[t^\frac12(\|H_t\|_{L^2}^2
+\|\nabla H\|_{L^2}^2)\big].
\end{align}
Meanwhile, we infer from \eqref{tt}, \eqref{dd1}, and \eqref{5.5} that
\begin{align}\label{k14}
&C\int_0^Tt^\frac12\big(\|\nabla H\|_{L^2}^4+\|\nabla u\|_{L^2}^4\|\nabla H\|_{L^2}^2+\|H\|_{L^2}\|\nabla H\|_{L^2}^3\big)dt\nonumber\\
&\le C\sup_{0\le t\le T}
\big(t^\frac12\|\nabla H\|_{L^2}+\|\nabla u\|_{L^2}^2\big)
\int_{0}^T\big(\|\nabla H\|_{L^2}^2+t^\frac12\|\nabla u\|_{L^2}^2\big)dt
\le C.
\end{align}
Substituting \eqref{k7}--\eqref{k14} into \eqref{k6} and choosing $\delta$ suitably small, we arrive at
\begin{align}\label{5.15}
&\sup_{0\le t\le T}\big[t^\frac12\big(\|\nabla u\|_{L^2}^2
+\|\nabla\theta\|_{L^2}^2+\|\sqrt{\rho}\dot{u}\|_{L^2}^2+\|H_t\|_{L^2}^2\big)\big]
\nonumber\\
&\quad+\int_0^Tt^\frac12\big(\|\sqrt{\rho}\dot{u}\|_{L^2}^2
+\|H_t\|_{L^2}^2+\|\nabla^2H\|_{L^2}^2+\|\sqrt{\rho}\dot{\theta}\|_{L^2}^2
+\|\nabla\dot{u}\|_{L^2}^2+\|\nabla H_t\|_{L^2}^2\big)dt\nonumber\\
&\le C,
\end{align}
provided that $$\mathbb{E}_0\le \varepsilon_0'\triangleq\min\left\{\varepsilon_4, \Big(\frac{(2c_2+\mu)\delta_1}{8C_{10}\mu}\Big)^\frac32, \Big(\frac{\mu}{12C_{11}}\Big)^\frac32, \Big(\frac{\eta}{6C_{11}}\Big)^\frac32\right\}.$$
Consequently, \eqref{1.11} follows from \eqref{5.5} and \eqref{5.15}.

\textbf{Step 4. Estimate for $\|H_t(t)\|_{L^2}$}. We infer from \eqref{rrf}, Sobolev's inequality, Gagliardo-Nirenberg inequality, Young's inequality, and \eqref{tt} that
\begin{align*}
& \frac12\frac{d}{dt}\|H_t\|_{L^2}^2+\eta \|\nabla H_t\|_{L^2}^2 \notag \\
&\le C\|\nabla u\|_{L^2}\|H_t\|_{L^4}^2
+C\|\nabla H_t\|_{L^2}\|H\|_{L^3}\|\dot{u}\|_{L^6}
+C\|\nabla H_t\|_{L^2}\|u\|_{L^6}\|H_t\|_{L^3}\nonumber\\
&\quad
+C\|H\|_{L^{12}}\|\nabla H_t\|_{L^2}\|u\|_{L^6}\|\nabla u\|_{L^4}\nonumber\\
&\le C\|\nabla u\|_{L^2}\|H_t\|_{L^2}^\frac12\|\nabla H_t\|_{L^2}^\frac32
+C\|\nabla H_t\|_{L^2}\|H\|_{L^2}^{\frac12}\|\nabla H\|_{L^2}^{\frac12}\|\nabla\dot{u}\|_{L^2}
\nonumber\\
&\quad
+C\|\nabla u\|_{L^2}\|\nabla u\|_{L^4}\|\nabla H\|_{L^2}^\frac34
|\nabla^2 H\|_{L^2}^\frac14\|\nabla H_t\|_{L^2}
\nonumber\\
&\le \frac{\eta}{2}\|\nabla H_t\|_{L^2}^2+C\|\nabla u\|_{L^2}^4\|H_t\|_{L^2}^2
+C\|\nabla\dot{u}\|_{L^2}^2\|\nabla H\|_{L^2}
\nonumber\\
&\quad+C\|\nabla u\|_{L^4}^4+C\|\nabla^2H\|_{L^2}^2\|\nabla H\|_{L^2}^2+C\|\nabla u\|_{L^2}^4\|\nabla H\|_{L^2}^2,
\end{align*}
which multiplied by $t$ leads to
\begin{align}\label{5.16}
\frac{d}{dt}\big(t\|H_t\|_{L^2}^2\big)+\eta t\|\nabla H_t\|_{L^2}^2
&\le Ct\|\nabla u\|_{L^2}^4\|H_t\|_{L^2}^2
+Ct\|\nabla\dot{u}\|_{L^2}^2\|\nabla H\|_{L^2}
+Ct\|\nabla^2H\|_{L^2}^2\|\nabla H\|_{L^2}^2\nonumber\\
&\quad +Ct\|\nabla u\|_{L^2}^4\|\nabla H\|_{L^2}^2
+Ct\|\nabla u\|_{L^4}^4+\|H_t\|_{L^2}^2.
\end{align}
We obtain from \eqref{a1}$_4$, Sobolev's inequality, and \eqref{tt} that
\begin{align}\label{5.18}
\|H_t\|_{L^2}^2 & \le C\|\nabla^2 H\|_{L^2}^2
+C\||u||\nabla H|\|_{L^2}^2+C\||\nabla u||H|\|_{L^2}^2 \nonumber\\
&\le C\|\nabla^2 H\|_{L^2}^2
+C\|u\|_{L^6}^2\|\nabla H\|_{L^3}^2+C\|H\|_{L^\infty}^2\|\nabla u\|_{L^2}^2 \nonumber\\
&\le C\|\nabla^2 H\|_{L^2}^2
+C\|\nabla u\|_{L^2}^2\|\nabla H\|_{H^1}^2+C\|H\|_{H^2}^2\|\nabla u\|_{L^2}^2
\nonumber\\
&\le C\|\nabla^2 H\|_{L^2}^2+C\|\nabla u\|_{L^2}^2.
\end{align}
Integrating \eqref{5.16} over $[0, T]$, we deduce from \eqref{5.17}, \eqref{5.18}, \eqref{dd1}, \eqref{5.5}, and \eqref{tt} that
\begin{align}\label{5.19}
&\sup_{0\le t\le T}\big(t\|H_t\|_{L^2}^2\big)+\int_0^Tt\|\nabla H_t\|_{L^2}^2dt\nonumber\\
&\le \int_0^T\|H_t\|_{L^2}^2dt+C\sup_{0\le t\le T}\big(t^{\frac12}\|\nabla u\|_{L^2}^2\big)^2\int_0^T\|H_t\|_{L^2}^2dt
+C\sup_{0\le t\le T}\big(t\|\nabla H\|_{L^2}^2\big)^\frac12\int_0^Tt^\frac12\|\nabla\dot{u}\|_{L^2}^2dt
\nonumber\\
&\quad
+C\sup_{0\le t\le T}\big(t\|\nabla H\|_{L^2}^2\big)\int_0^T\|\nabla^2H\|_{L^2}^2dt
+C\sup_{0\le t\le T}\big(t^{\frac12}\|\nabla u\|_{L^2}^2\big)^2\int_0^T\|\nabla H\|_{L^2}^2dt
\nonumber\\
&\quad+\sup_{0\le t\le T}\big[t^\frac12\big(\|\nabla u\|_{L^2}^2+\|\nabla\theta\|_{L^2}^2+\|\sqrt{\rho}\dot{u}\|_{L^2}^2\big)\big]
\int_0^Tt^\frac12\big(\|\nabla u\|_{L^2}^2
+\|\nabla\theta\|_{L^2}^2+\|\sqrt{\rho}\dot{u}\|_{L^2}^2\big)dt\nonumber\\
&\le C.
\end{align}

\textbf{Step 5. Estimate for $\|\sqrt{\rho}\dot{\theta}(t)\|_{L^2}$}. We deduce from \eqref{f80}, H\"older's inequality, Sobolev's inequality, \eqref{3.49}, \eqref{f47}, \eqref{5.17}, \eqref{3.74}, and \eqref{tt} that
\begin{align*}
& c_v\frac{d}{dt}\|\sqrt{\rho}\dot{\theta}\|_{L^2}^2+\kappa
\|\nabla\dot{\theta}\|_{L^2}^2\nonumber\\
&\le C\big(\|\sqrt{\rho}\theta\|_{L^2}\|\sqrt{\rho}\theta\|_{L^6}
+\|\nabla u\|_{L^2}\|\nabla u\|_{L^6}+\|\nabla H\|_{H^1}^2\big)\nonumber\\
&\quad \times
\big(\|\sqrt{\rho}\dot{\theta}\|_{L^2}^2+\|\nabla^2\theta\|_{L^2}^2+\|\nabla H_t\|_{L^2}^2+\|H_t\|_{L^2}^2
+\|\nabla u\|_{L^4}^4+\|\nabla\dot{u}\|_{L^2}^2\big)\nonumber\\
&\le C\big(\|\nabla\theta\|_{L^2}^2+\|\nabla u\|_{L^2}^2
+\|\sqrt{\rho}\dot{u}\|_{L^2}^2+\|\nabla H\|_{L^2}^2
+\|\nabla^2 H\|_{L^2}^2\big)\nonumber\\
&\quad\times\big(\|\sqrt{\rho}\dot{\theta}\|_{L^2}^2
+\|\nabla\theta\|_{L^2}^2+\|\nabla u\|_{L^2}^2
+\|\sqrt{\rho}\dot{u}\|_{L^2}^2+\|\nabla H\|_{L^2}^2+\|\nabla H_t\|_{L^2}^2+\|H_t\|_{L^2}^2+\|\nabla\dot{u}\|_{L^2}^2\big),
\end{align*}
which multiplied by $t$ gives that
\begin{align}\label{5.20}
& c_v\frac{d}{dt}\big(t\|\sqrt{\rho}\dot{\theta}\|_{L^2}^2\big)
+\kappa t\|\nabla\dot{\theta}\|_{L^2}^2\nonumber\\
&\le C\big[t^{\frac12}\big(\|\nabla\theta\|_{L^2}^2+\|\nabla u\|_{L^2}^2
+\|\sqrt{\rho}\dot{u}\|_{L^2}^2\big)\big]\nonumber\\
&\quad\times t^{\frac12}\big(\|\sqrt{\rho}\dot{\theta}\|_{L^2}^2
+\|\nabla\theta\|_{L^2}^2+\|\nabla u\|_{L^2}^2
+\|\sqrt{\rho}\dot{u}\|_{L^2}^2+\|\nabla H_t\|_{L^2}^2+\|H_t\|_{L^2}^2+\|\nabla\dot{u}\|_{L^2}^2\big) \notag \\
& \quad +C\big(t\|\nabla H\|_{L^2}^2\big)
\big(\|\nabla\theta\|_{L^2}^2+\|\nabla u\|_{L^2}^2+\|\nabla\dot{\theta}\|_{L^2}^2
+\|\nabla\dot{u}\|_{L^2}^2+\|\nabla H\|_{H^1}^2+\|\nabla H_t\|_{H^1}^2\big)
 \notag \\
& \quad +C\|\nabla H\|_{H^1}^2
\big[t\big(\|H_t\|_{L^2}^2+\|\nabla H_t\|_{L^2}^2\big)\big]
+C\|\nabla\dot{\theta}\|_{L^2}^2 \notag \\
& \quad +C\big(\|\sqrt{\rho}\dot{\theta}\|_{L^2}^2
+\|\nabla\theta\|_{L^2}^2+\|\nabla u\|_{L^2}^2
+\|\sqrt{\rho}\dot{u}\|_{L^2}^2\big)\big(t\|\nabla^2 H\|_{L^2}^2\big) \notag \\
& \quad
+C\big[t\big(\|H_t\|_{L^2}^2+\|\nabla H\|_{L^2}^2\big)\big]\|\nabla\dot{u}\|_{L^2}^2
\end{align}
where we have used \eqref{f49} and the following facts
\begin{align*}
\|\sqrt{\rho}\dot{u}\|_{L^2}^2
\leq \|\sqrt{\rho}\|_{L^3}^2\|\dot{u}\|_{L^6}^2
\leq C(\|\rho_0\|_{L^1},\bar\rho)\|\nabla\dot{u}\|_{L^2}^2,\
\|\sqrt{\rho}\dot{\theta}\|_{L^2}^2
\leq \|\sqrt{\rho}\|_{L^3}^2\|\dot{\theta}\|_{L^6}^2
\leq C(\|\rho_0\|_{L^1},\bar\rho)\|\nabla\dot{\theta}\|_{L^2}^2,
\end{align*}
due to H\"older's inequality, \eqref{lz3.1}, \eqref{3.49}, and Sobolev's inequality.
As a consequence, integrating \eqref{5.20} over $[0,T]$, we deduce from \eqref{5.15}, \eqref{5.5}, \eqref{dd1}, \eqref{tt}, and \eqref{5.19} that
\begin{align}\label{5.21}
\sup_{0\le t\le T}\big(t\|\sqrt{\rho}\dot{\theta}\|_{L^2}^2\big)
+\int_0^Tt\|\nabla\dot{\theta}\|_{L^2}^2dt\le C.
\end{align}

\textbf{Step 6. Estimate for $\|\nabla^2H(t)\|_{L^2}$ and $\|\nabla^2\theta(t)\|_{L^2}$}. It follows from \eqref{f49}, \eqref{5.19}, \eqref{dd1}, and \eqref{5.15} that
\begin{align}\label{5.22}
\sup_{0\le t\le T}\big(t\|\nabla^2H\|_{L^2}^2\big)
\le C\sup_{0\le t\le T}\big(t\|H_t\|_{L^2}^2\big)
+C\sup_{0\le t\le T}\big(t^{\frac12}\|\nabla H\|_{L^2}\big)\sup_{0\le t\le T}\big(t^{\frac12}\|\nabla u\|_{L^2}^2\big)\le C.
\end{align}
We obtain from \eqref{a1}$_3$, the standard $L^2$-estimates of elliptic equations, Sobolev's inequality, Gagliardo-Nirenberg inequality, \eqref{tt}, and \eqref{f50} that
\begin{align*}
\|\nabla^2\theta\|_{L^2}^2&\le C\big(\|\sqrt{\rho}\dot{\theta}\|_{L^2}^2+\|\theta\nabla u\|_{L^2}^2+\|\nabla u\|_{L^4}^4
+\|\nabla H\|_{L^4}^4\big)\nonumber\\
&\le C\big(\|\sqrt{\rho}\dot{\theta}\|_{L^2}^2+\|\theta\|_{L^\infty}^2\|\nabla u\|_{L^2}^2+\|\nabla u\|_{L^2}\|\nabla u\|_{L^6}^3
+\|\nabla H\|_{L^2}\|\nabla^2H\|_{L^2}^3\big)\nonumber\\
&\le C\|\sqrt{\rho}\dot{\theta}\|_{L^2}^2+C\|\nabla\theta\|_{L^2}\|\nabla^2\theta\|_{L^2}\|\nabla u\|_{L^2}^2
+C\|\nabla^2H\|_{L^2}^2\nonumber\\
&\quad
+C\|\nabla u\|_{L^2}\big(\|\sqrt{\rho}\dot{u}\|_{L^2}^3
+\|H_t\|_{L^2}^3+\|\nabla\theta\|_{L^2}^3+\|\nabla H\|_{L^2}^3\big)
\nonumber\\
&\le C\|\sqrt{\rho}\dot{\theta}\|_{L^2}^2
+\frac12\|\nabla^2\theta\|_{L^2}^2
+C\|\nabla\theta\|_{L^2}^2\|\nabla u\|_{L^2}^4+C\|\nabla^2H\|_{L^2}^2
\nonumber\\ &\quad +C\big(\|\nabla u\|_{L^2}^4+\|\sqrt{\rho}\dot{u}\|_{L^2}^4
+\|H_t\|_{L^2}^4+\|\nabla\theta\|_{L^2}^4+\|\nabla H\|_{L^2}^4\big).
\end{align*}
This along with \eqref{5.21}, \eqref{5.22}, and \eqref{5.15} implies that
\begin{align*}
\sup_{0\le t\le T}\big(t\|\nabla^2\theta\|_{L^2}^2\big)
&\le C\sup_{0\le t\le T}\big(t\|\sqrt{\rho}\dot{\theta}\|_{L^2}^2\big)
+C\sup_{0\le t\le T}\big(t\|\nabla^2H\|_{L^2}^2\big) \nonumber\\
&\quad
+C\sup_{0\le t\le T}\big[t^\frac12\big(\|\nabla u\|_{L^2}^2
+\|\nabla\theta\|_{L^2}^2+\|\sqrt{\rho}\dot{u}\|_{L^2}^2
+\|H_t\|_{L^2}^2\big)\big]^2\le C,
\end{align*}
which combined with \eqref{dd1}, \eqref{5.19}, \eqref{5.21}, and \eqref{5.22} yields \eqref{1.12}.
\hfill $\Box$

\section*{Conflict of interest}
The authors declare that they have no conflict of interest.

\section*{Data availability}
Data sharing is not applicable to this article as no new data were created or analyzed in this study.


\begin{thebibliography}{10}


\renewcommand\refname{References}
\renewenvironment{thebibliography}[1]{%
\section*{\refname}
\list{{\arabic{enumi}}}{\def\makelabel##1{\hss{##1}}\topsep=0mm
\parsep=0mm
\partopsep=0mm\itemsep=0mm
\labelsep=1ex\itemindent=0mm
\settowidth\labelwidth{\small[#1]}%
\leftmargin\labelwidth \advance\leftmargin\labelsep
\advance\leftmargin -\itemindent
\usecounter{enumi}}\small
\def\newblock{\ }
\sloppy\clubpenalty4000\widowpenalty4000
\sfcode`\.=1000\relax}{\endlist}
\bibliographystyle{model1b-num-names}

\bibitem{D2017}
P. A. Davidson, Introduction to magnetohydrodynamics, 2nd edition, Cambridge University Press, Cambridge, 2017.

\bibitem{D1997}
B. Desjardins, Regularity of weak solutions of the compressible isentropic Navier-Stokes equations, {\it Comm. Partial Differential Equations}, {\bf 22} (1997), 977--1008.

\bibitem{D2016}
E. DiBenedetto, Real analysis, 2nd edition, Birkh{\"a}user, New York, 2016.

\bibitem{DF2006}
B. Ducomet and E. Feireisl, The equations of magnetohydrodynamics: on the interaction between matter and radiation in the evolution of gaseous stars, {\it Comm. Math. Phys.}, {\bf 266} (2006), 595--629.

\bibitem{FY09}
J. Fan and W. Yu, Strong solution to the compressible magnetohydrodynamic equations with vacuum, {\it Nonlinear Anal. Real World Appl.}, {\bf 10} (2009), 392--409.

\bibitem{F2004}
E. Feireisl, Dynamics of viscous compressible fluids, Oxford University Press, Oxford, 2004.

\bibitem{FL2020}
E. Feireisl and Y. Li, On global-in-time weak solutions to the magnetohydrodynamic system of compressible inviscid fluids, {\it Nonlinearity}, {\bf 33} (2020), 139--155.

\bibitem{HHPZ}
G. Hong, X. Hou, H. Peng, and C. Zhu, Global existence for a class of large solutions to three dimensional compressible magnetohydrodynamic equations with vacuum, {\it SIAM J. Math. Anal.}, {\bf 49} (2017), 2409--2441.

\bibitem{HJP22}
X. Hou, M. Jiang, and H. Peng, Global strong solution to 3D full compressible magnetohydrodynamic flows with vacuum at infinity, {\it Z. Angew. Math. Phys.}, {\bf 73} (2022), Paper No. 13.

\bibitem{HW08}
X. Hu and D. Wang, Global solutions to the three-dimensional full compressible magnetohydrodynamic flows, {\it Commun. Math. Phys.}, {\bf 283} (2008), 255--284.

\bibitem{HW10}
X. Hu and D. Wang, Global existence and large-time behavior of solutions to the three dimensional equations of compressible magnetohydrodynamic flows, {\it Arch. Ration. Mech. Anal.}, {\bf 197} (2010), 203--238.

\bibitem{HL13}
X. Huang and J. Li, Serrin-type blowup criterion for viscous, compressible, and heat conducting Navier-Stokes and magnetohydrodynamic flows, {\it Comm. Math. Phys.}, {\bf 324} (2013), 147--171.

\bibitem{HL18}
X. Huang and J. Li, Global classical and weak solutions to the three-dimensional full compressible Navier-Stokes system with vacuum and large oscillations, {\it Arch. Ration. Mech. Anal.}, {\bf 227} (2018), 995--1059.

\bibitem{HLX20112}
X. Huang, J. Li, and Z. Xin, Serrin-type criterion for the
three-dimensional viscous compressible flows, {\it SIAM J. Math. Anal.},
{\bf 43} (2011), 1872--1886.

\bibitem{JJL12}
S. Jiang, Q. Ju, and Y. Li, Low Mach number limit for the multi-dimensional full magnetohydrodynamic equations, {\it Nonlinearity}, {\bf 25} (2012), 1351--1365.

\bibitem{JJLX14}
S. Jiang, Q. Ju, F. Li, and Z. Xin, Low Mach number limit for the full compressible magnetohydrodynamic equations with general initial data, {\it Adv. Math.}, {\bf 259} (2014), 384--420.

\bibitem{LXZ13}
H. Li, X. Xu, and J. Zhang, Global classical solutions to 3D compressible magnetohydrodynamic equations with large oscillations and vacuum, {\it SIAM J. Math. Anal.}, {\bf 45} (2013), 1356--1387.

\bibitem{JL19}
J. Li, Global small solutions of heat conductive compressible Navier-Stokes equations with vaccum: smallness on scaling invariant quantity,
{\it Arch. Ration. Mech. Anal.}, {\bf 237} (2020), 899--919.

\bibitem{LQ2012}
T. Li and T. Qin, Physics and partial differential equations. vol. 1, Translated from the Chinese original by Yachun Li, Higher Education Press, Beijing, 2012.

\bibitem{LG2014}
X. Li and B. Guo, On the equations of thermally radiative magnetohydrodynamics, {\it J. Differential Equations}, {\bf 257} (2014), 3334--3381.

\bibitem{LS2021}
Y. Li and Y. Sun, On global-in-time weak solutions to a two-dimensional full compressible non-resistive MHD system, {\it SIAM J. Math. Anal.}, {\bf 53} (2021), 4142--4177.

\bibitem{L21}
Z. Liang, Global strong solutions of Navier-Stokes equations for heat-conducting compressible fluids with vacuum at infinity,
{\it J. Math. Fluid Mech.}, {\bf 23} (2021), Paper No. 17.

\bibitem{L1998}
P. L. Lions, Mathematical topics in fluid mechanics, vol. II: compressible models, Oxford University Press, Oxford, 1998.

\bibitem{LZ20}
Y. Liu and X. Zhong, Global well-posedness to three-dimensional full compressible magnetohydrodynamic equations with vacuum, {\it Z. Angew. Math. Phys.}, {\bf 71} (2020), Paper No. 188.

\bibitem{LZ22}
Y. Liu and X. Zhong, Global strong solution for 3D compressible heat-conducting magnetohydrodynamic equations revisited, {\it J. Differential Equations}, {\bf 336} (2022), 456--478.

\bibitem{LSX16}
B. L{\"u}, X. Shi, and X. Xu, Global existence and large-time asymptotic behavior of strong solutions to the compressible magnetohydrodynamic equations with vacuum, {\it Indiana Univ. Math. J.}, {\bf 65} (2016), 925--975.

\bibitem{S70}
E. M. Stein, Singular integrals and differentiability properties of functions, Princeton University Press, Princeton, N.J., 1970.

\bibitem{SH12}
A. Suen and D. Hoff, Global low-energy weak solutions of the equations of three-dimensional compressible magnetohydrodynamics, {\it Arch. Ration. Mech. Anal.}, {\bf 205} (2012), 27--58.

\bibitem{T1976}
G. Talenti, Best constants in Sobolev inequality, {\it Ann. Mat. Pura Appl.},
{\bf 110} (1976), 353--372.

\bibitem{W21}
Y. Wang, A Beale-Kato-Majda criterion for three dimensional compressible viscous non-isentropic magnetohydrodynamic flows without heat-conductivity, {\it J. Differential Equations}, {\bf 280} (2021), 66--98.

\bibitem{WZ17}
H. Wen and C. Zhu, Global solutions to the three-dimensional full compressible Navier-Stokes equations with vacuum at infinity in some classes of large data, {\it SIAM J. Math. Anal.}, {\bf 49} (2017), 162--221.

\bibitem{WW17}
J. Wu and Y. Wu, Global small solutions to the compressible 2D magnetohydrodynamic system without magnetic diffusion, {\it Adv. Math.}, {\bf 310} (2017), 759--888.

\bibitem{XZ21}
H. Xu and J. Zhang, Regularity and uniqueness for the compressible full Navier-Stokes equations, {\it J. Differential Equations}, {\bf 272} (2021), 46--73.

\bibitem{YZ20}
H. Yu and P. Zhang, Global strong solutions to the 3D full compressible Navier-Stokes equations with density-temperature-dependent viscosities in bounded domains, {\it J. Differential Equations}, {\bf 268} (2020), 7286--7310.

\bibitem{Z19}
X. Zhong, On formation of singularity of the full compressible magnetohydrodynamic equations with zero heat conduction, {\it Indiana Univ. Math. J.}, {\bf 68} (2019), 1379--1407.

\end{thebibliography}
\end{document}